\documentclass[11pt]{article}
\usepackage{amsmath,amssymb,amsthm,mathtools}
\ifdefined\acmJournal\else
\usepackage[margin=1in]{geometry}
\usepackage[numbers]{natbib}
\usepackage[hidelinks]{hyperref}
\usepackage[expansion=false]{microtype}
\fi
\newtheorem{theorem}{Theorem}[section]
\newtheorem{lemma}[theorem]{Lemma}
\newtheorem{proposition}[theorem]{Proposition}

\numberwithin{equation}{section}
\makeatletter
\@ifclassloaded{article}{%
  \newcommand{\acmJournal}[1]{}
  \newcommand{\orcid}[1]{\def\preprintORCID{#1}}
  \newcommand{\email}[1]{\def\preprintEmail{#1}}
  \newcommand{\department}[1]{#1, }
  \newcommand{\institution}[1]{#1, }
  \newcommand{\city}[1]{#1 }
  \newcommand{\postcode}[1]{#1, }
  \newcommand{\country}[1]{#1}
  \newcommand{\affiliation}[1]{\g@addto@macro\@author{\thanks{Corresponding author.
    #1. Email: \href{mailto:\preprintEmail}{\texttt{\preprintEmail}}.
    \href{https://orcid.org/\preprintORCID}{ORCID: \preprintORCID}.}}}
  \NewDocumentEnvironment{CCSXML}{+b}{}{}
  \let\preprintConcepts\@empty
  \def\preprintConcept#1~#2\@nil{#1 \(\rightarrow\) #2}
  \newcommand{\ccsdesc}[2][]{%
    \ifx\preprintConcepts\@empty\else\g@addto@macro\preprintConcepts{; }\fi
    \g@addto@macro\preprintConcepts{\preprintConcept#2\@nil}}
  \newcommand{\keywords}[1]{\def\preprintKeywords{#1}\hypersetup{pdfkeywords={#1}}}
  \AddToHook{env/thebibliography/begin}{\raggedright}

  \RenewDocumentEnvironment{abstract}{+b}{\gdef\preprintAbstract{#1}}{}
  \AddToHook{cmd/maketitle/after}{%
    \begin{articleabstract}
      \preprintAbstract\par\smallskip
      \noindent\textbf{CCS Concepts.} \preprintConcepts.\par
      \noindent\textbf{Keywords.} \preprintKeywords.
    \end{articleabstract}}
}{}
\makeatother

\acmJournal{JACM}
\title{Deterministic regularisation of the nonsymmetric eigenproblem}
\author{Matthew J. Colbrook}
\orcid{0000-0003-4964-9575}
\affiliation{%
  \department{Department of Applied Mathematics and Theoretical Physics}%
  \institution{University of Cambridge}%
  \city{Cambridge}%
  \postcode{CB3 0WA}%
  \country{United Kingdom}}
\email{m.colbrook@damtp.cam.ac.uk}
\date{}
\hypersetup{pdftitle={Deterministic regularisation of the nonsymmetric eigenproblem},
pdfauthor={Matthew J. Colbrook}}
\begin{document}

\begin{abstract}
We give a deterministic algorithm for pseudospectral shattering in exact arithmetic, resolving an important open problem in the complexity of the nonsymmetric eigenproblem. Given a complex matrix $A$ of order $n\ge2$ with $\|A\|\le1$
and a tolerance $0<\delta<1/2$, the algorithm computes a perturbation
$E$ of norm less than $\delta$ and an exact eigensystem of $A+E$ in
$O(n^3\operatorname{polylog}(n/\delta))$ real-field operations and
comparisons. The spectrum of $A+E$ is simple, and the ratio of its
eigenbasis condition number to its minimum eigenvalue gap is at most
$C(n/\delta)^{20/3}$ for an absolute constant $C$.
To choose the perturbation, we study how the eigenvalues of a rank-one
pencil $B+tuv^*$ move as $t$ varies. If the initial derivative of an
eigenvalue is large in modulus compared with its distance to the rest
of the spectrum, then two eigenvalues coincide at a complex value of $t$
near zero. The discriminant vanishes at each such collision. Jensen's
formula then gives a lower bound for its mean logarithmic modulus on a
circle about zero. We use this estimate to choose successive
perturbations while controlling the norm of their sum.
\end{abstract}

\ccsdesc[500]{Mathematics of computing~Computations on matrices}
\ccsdesc[300]{Computing methodologies~Linear algebra algorithms}
\ccsdesc[300]{Theory of computation~Pseudorandomness and derandomization}
\keywords{pseudospectral shattering, nonsymmetric eigenproblem, deterministic algorithms, spectral perturbation theory}
\maketitle

\section{Introduction}

A small perturbation of a matrix can separate a multiple eigenvalue
and still leave the eigenvectors nearly parallel. For example, consider
\[
 B_t=\begin{pmatrix}0&1\\ t&0\end{pmatrix},\qquad 0<t<1.
\]
Its eigenvalues are $\pm\sqrt t$, and the corresponding eigenvectors
$(1,\pm\sqrt t)^T$ approach the same vector as $t\downarrow0$.
A diagonal matrix, on the other hand, has an orthonormal eigenbasis
irrespective of how close its eigenvalues are. To regularise the
nonsymmetric eigenproblem we therefore seek a nearby matrix for which
we can bound both the eigenbasis condition number and the reciprocal
of the minimum eigenvalue gap.

For a matrix of order $n$ and norm at most one, a Gaussian perturbation
of size $\delta$ gives polynomial bounds in $n/\delta$ on both
quantities with high probability~\cite{gaussian,banks}.
Can we find such a perturbation deterministically in cubic time up to
logarithmic factors? The algorithm of Theorem~\ref{thm:main} does so in
exact arithmetic. It also computes an exact eigensystem of the perturbed
matrix, even when the input is defective.

Davies~\cite{davies} obtained bounds for approximate diagonalisation
with exponents that depend on the dimension. Banks, Kulkarni,
Mukherjee and Srivastava~\cite{gaussian} proved polynomial bounds on
the eigenbasis condition number under Gaussian perturbations.
Banks, Garza-Vargas, Kulkarni and Srivastava~\cite{banks} also bounded
the reciprocal gap. They used the resulting separation of the
pseudospectrum to obtain a numerically stable randomised diagonalisation
algorithm in nearly matrix multiplication time. For matrices of
polynomially bounded norm, Shah, Srivastava and Zeng~\cite[Theorem~1.1]{shah}
obtained polynomial bounds on both quantities with high probability
for Bernoulli--Gaussian perturbations of expected support size
$n^{1+\alpha}$, where $0<\alpha\le1/2$ is fixed.

Banks et al.~\cite[Section~6]{banks} raised the question of a
deterministic construction. Amsel et al.~\cite[Problem~3.1 and
equation~(6)]{workshop} gave the quantitative formulation. The problem
also appears as Problem~IE-07 in Townsend's collection~\cite{townsend}.

\subsection{Main result and computational model}\label{subsec:main-model}

Let $\|\cdot\|$ denote the Euclidean norm and its induced matrix norm,
and let $\|\cdot\|_F$ denote the Frobenius norm. For an invertible
matrix $X$, its condition number is $\kappa_2(X)=\|X\|\,\|X^{-1}\|$.
If $B$ has distinct eigenvalues $\lambda_1,\ldots,\lambda_n$, with
$n\ge2$, define
\[
 \operatorname{gap}(B)=\min_{i\ne j}|\lambda_i-\lambda_j|,
 \qquad
 \kappa_V(B)=\inf_{B=X\Lambda X^{-1}}\kappa_2(X),
\]
where the infimum is over all diagonal $\Lambda$ and invertible $X$.
For a matrix of order one, $\kappa_V(B)=1$ and no minimum gap is needed.

\begin{theorem}[Deterministic spectral regularisation]\label{thm:main}
Let $A\in\mathbb C^{n\times n}$ satisfy $\|A\|\le1$, where $n\ge2$,
and let $0<\delta<1/2$. A deterministic algorithm computes a
perturbation $E$ with $\|E\|<\delta$ such that $B=A+E$ has simple
spectrum and
\[
  \frac{\kappa_V(B)}{\operatorname{gap}(B)}
     \le C\left(\frac{n}{\delta}\right)^{20/3}.
\]
The algorithm also computes the eigenvalues of $B$, an eigenbasis and
its inverse. The total cost is $O(n^3\log^{c_{\log}}(n/\delta))$ real-field
operations and comparisons, where $C$ and $c_{\log}$ are absolute constants.
\end{theorem}

To see what the bound implies for the pseudospectrum, let
$\operatorname{spec}(B)$ denote the spectrum of $B$. For $\varepsilon>0$,
the $\varepsilon$-pseudospectrum is
\[
 \operatorname{spec}(B)\,\cup\,
 \{z\notin\operatorname{spec}(B):\|(zI-B)^{-1}\|>\varepsilon^{-1}\}.
\]
Diagonalisation gives
\[
 \|(zI-B)^{-1}\|\le
 \frac{\kappa_V(B)}{\min_i|z-\lambda_i|}
 \qquad(z\notin\operatorname{spec}(B)).
\]
Hence the $\varepsilon$-pseudospectrum is contained in the union of the
discs of radius $\varepsilon\kappa_V(B)$ about the eigenvalues. If
$\varepsilon<\operatorname{gap}(B)/(2\kappa_V(B))$, these discs are
pairwise disjoint. For the matrix in Theorem~\ref{thm:main}, this holds
whenever $0<\varepsilon<(2C)^{-1}(\delta/n)^{20/3}$. Each disc then
contains one eigenvalue of $B+F$ for every further perturbation
$\|F\|<\varepsilon$, by continuity of the eigenvalues along
$B+tF$, $0\le t\le1$. For background on pseudospectra,
see Trefethen and Embree~\cite{trefethen-embree}.

The formulation in~\cite{workshop} uses an eigenbasis with unit columns.
If $V$ is any eigenbasis and $D$ is the diagonal matrix of its column
norms, then $\|VD^{-1}\|\le\sqrt n$ and
$\|DV^{-1}\|\le\|V\|\,\|V^{-1}\|$. Thus normalising the columns costs
at most a factor of $\sqrt n$ in the condition number bound.

\paragraph{Computational model.}
Each real addition, subtraction, multiplication, division by a nonzero
scalar and comparison has unit cost. A complex scalar is represented by
two real scalars. In particular, the cost of arithmetic on integers is
independent of their bit length. Section~\ref{sec:mixed-range-bridge}
uses this convention to encode a column of integers as the digits of
one integer. Binary searches over integers use intervals of dyadic
length, extending the original interval if necessary.
The algorithm approximates polynomial roots by field operations.
A further perturbation makes these approximations exact eigenvalues
and the corresponding computed vectors exact eigenvectors. Its norm
is included in the perturbation bound of Theorem~\ref{thm:main}.

Here and below $M^*$ is the conjugate transpose of $M$,
$\chi_M(z)=\det(zI-M)$ is its characteristic polynomial, and $I$ is
the identity of the indicated order. The closed disc of centre $c$ and
radius $r$ is denoted by $D(c,r)$. A contraction is a matrix of norm
at most one. All logarithms without a subscript are natural. A fixed
ordering of the indices resolves ties in finite searches.

\subsection{Proof strategy}
We start from a nearby matrix $B$ with distinct eigenvalues and a known
eigenbasis. Let $P_i$ be the spectral projector for $\lambda_i$ and
$s_i$ the distance from $\lambda_i$ to the rest of the spectrum.
The norm $\|P_i\|$ measures the first-order sensitivity of $\lambda_i$.
For the pencil $B+tuv^*$, the derivative of this eigenvalue at $t=0$
is $v^*P_i u$. The first step is to choose $u$ and $v$
deterministically so that $|v^*P_i u|$ is large relative to $s_i$
for sufficiently many of the indices with large $\|P_i\|/s_i$.

For this choice of directions, the discriminant of the pencil has zeros
near $t=0$. Jensen's formula expresses the mean of the logarithm of
its modulus on a circle in terms of these zeros. We subtract a multiple
of the squared Frobenius norm of the accumulated perturbation, using
an approximation to that perturbation when we reuse an eigenbasis.
The term in this penalty that is linear in $t$ has mean zero on the
circle. The mean estimate therefore gives a choice of $t$ that
increases the potential while keeping the matrix close to the input.

Recomputing a full eigenbasis after every rank-one perturbation would
be too expensive. The decision to recompute it depends on a sum of
logarithmic weights derived from the ratios $\|P_i\|/s_i$. When the
sum is large, the lower bound on the increase in the potential limits
the number of times a new basis is needed. When the sum is small, it
bounds the total precision needed for the reference eigenvectors.
Several updates can then use one basis. We estimate all the projector
norms in the original Euclidean norm but recover left and right
eigenvectors only for the indices used to choose the next perturbation.

When only a few ratios exceed the prescribed bound, the remaining
right eigenvectors span an invariant subspace of small codimension.
In a basis adapted to that subspace, the matrix is block upper
triangular. We apply the algorithm recursively to the trailing block.
A scalar shift separates the spectrum of this block from each retained
eigenvalue by an amount that depends on its original gap $s_i$.
This separation bounds the projector norms of the assembled matrix.
Two levels of recursion suffice for the operation count in
Theorem~\ref{thm:main}.

\subsection{Outline}
Section~\ref{sec:foundations} relates the growth of the discriminant
to the choice of $u$ and $v$. Section~\ref{sec:scalar-batch} gives the
polynomial algorithms for choosing the perturbation parameter and
approximating the roots of the characteristic polynomial.
Section~\ref{sec:seed} constructs the initial eigensystem.
Sections~\ref{sec:transfer-rebase} and~\ref{sec:mixed-range-bridge}
show how to recover eigenvectors and estimate projector norms over a
sequence of rank-one updates. Section~\ref{sec:global-assembly}
combines these results to prove Theorem~\ref{thm:main}.

\subsection{AI declaration}

During the development of this work, the author used ChatGPT 6 (OpenAI)
as an interactive aid. AI assistance was used to explore the resolvent-vector accuracy estimates in Lemma 5.3, the exact-diagonalisation correction and structured inversion in Proposition 5.4, the triangular factorisation with individual accuracies in Lemma 6.4, and the norm-estimation and direction-recovery constructions in Propositions 6.5 and 6.6. The author
specified the questions, hypotheses, and estimates required in these
exchanges, and used the suggestions as material for further mathematical
development. The author determined the final statements and
conclusions and takes full responsibility for  all content.

\section{Discriminant growth under a rank-one perturbation}
\label{sec:foundations}

For each simple eigenvalue, we compare the norm of its spectral projector
with the distance to the rest of the spectrum.
Let $B\in\mathbb C^{n\times n}$, $n\ge2$, have simple spectrum,
let $P_i$ denote the spectral projector for $\lambda_i$, and define
\[
 s_i=\min_{j\ne i}|\lambda_i-\lambda_j|,\qquad
 K_i=\frac{\|P_i\|}{s_i},\qquad
 D(B)=\prod_{i<j}(\lambda_i-\lambda_j)^2 .
\]
For a monic polynomial $p(z)=\prod_i(z-\zeta_i)$, we use the convention
$\operatorname{disc}(p)=\prod_{i<j}(\zeta_i-\zeta_j)^2$.
Thus $D(B)=\operatorname{disc}(\chi_B)$.
There are $d=n(n-1)$ linear factors in this product, counting
multiplicities. The formula also defines the discriminant for a repeated
spectrum, in which case it vanishes. For $x\ge0$, define $\log_2^+x=\log_2\max\{1,x\}$.
If $B=X\Lambda X^{-1}$ and the columns $x_i$ of $X$ have norm one,
the dual rows $w_i^*$ have norm $\|P_i\|$, since $P_i=x_iw_i^*$.
Bounding the norms of the eigenbasis and its inverse by their Frobenius norms gives
\begin{equation}\label{eq:basis-projector}
  \kappa_V(B)\le n\max_i\|P_i\|.
\end{equation}
Each projector satisfies $\|P_i\|\ge1$, and $s_i<4$ if $\|B\|<2$.

Approximate normalisation is sufficient for the algorithm.
For $0<\epsilon\le1$, we can scale a nonzero column to have norm in
$[(1+\epsilon)^{-1},1]$ as follows. Dividing by the largest absolute
real or imaginary component gives a squared norm in $[1,2n]$.
Bisection, using comparisons of squared lengths, then gives an upper
estimate of the length within relative error $\epsilon$. The total cost
is $O(n+\log(n/\epsilon))$. Scaling the dual row by the reciprocal
factor preserves the projector exactly.

\begin{lemma}[Collisions of a rank-one pencil]\label{lem:collision}
Suppose $B$ has simple spectrum. Let $B(t)=B+tuv^*$ and
$a_i=|v^*P_i u|$.
For every $r>0$, the polynomial $g(t)=D(B(t))$ has at least
$\tfrac12\#\{i:ra_i>4s_i\}$ zeros in $|t|\le2r$, counted with
multiplicity. Consequently, for every $h_0>0$,
\begin{equation}
 \log_2\frac{\max_{|t|=h_0}|g(t)|}{|g(0)|}
 \ \ge\ \frac{1}{2\pi}\int_0^{2\pi}
       \log_2\frac{|g(h_0e^{i\theta})|}{|g(0)|}\,d\theta
 \ge\frac12\sum_i\log_2^+\frac{h_0a_i}{8s_i}.
 \label{eq:ascent}
\end{equation}
\end{lemma}
\begin{proof}
By the determinant lemma, the characteristic polynomial of $B(t)$ is
$p(z)-tq(z)$, where
\[
 p(z)=\chi_B(z),\qquad q(z)=v^*\operatorname{adj}(zI-B)u.
\]
For the given $r$, call an index $i$ marked if $ra_i>4s_i$. If $q=0$,
then $g$ is constant and every $v^*P_i u$ vanishes, so both conclusions
follow. Otherwise, let
$h=\gcd(p,q)$, $f=p/h$ and $v_0=q/h$. Each root of $h$ gives a stationary eigenvalue branch, whose derivative is zero. Choose $R\in(r,2r)$ so that $|t|=R$ contains no critical value of
$f/v_0$. Since this rational map takes infinity to infinity, each
component of the preimage of $|t|<R$ is a plane domain and a proper
branched cover of the disc of some degree $j$. Each component of degree $j$ contains exactly $j$ roots
of $f$, since zero is a regular value. In a component of degree
$j\ge2$, the ramification is at least $j-1$. Indeed, deleting the
branch values gives a connected covering, whose monodromy is
generated by loops around those values and is transitive on the $j$ sheets. A cycle of length $e$ can join at most $e-1$ previously
disjoint classes of sheets. At a point of local degree $e$, the map has the form
$\zeta\mapsto\zeta^e$ in local coordinates. Its $e$ inverse branches
have pairwise differences of order $1/e$, so their squared differences
contribute order $2\binom e2/e=e-1$ to the discriminant. Since $j-1\ge j/2$, the ramification in these
components is at least half the number of marked roots they contain.

For a component of degree one, the inverse branch through
$\lambda_i$ is univalent. Differentiating the simple eigenvalue equation
gives its derivative $v^*P_i u$ at zero. Koebe's quarter
theorem shows that the open disc of radius $Ra_i/4$ about $\lambda_i$ lies in the
image of the inverse branch \cite[Chapter~2]{koebe}. For a marked index $i$, this disc
contains a nearest other eigenvalue of $B$. No second preimage of
zero can belong to the component, so this other eigenvalue is a root
of $h$. At its image under $f/v_0$ the moving and stationary
branches meet. This gives a zero of $\operatorname{Res}(h,f-tv_0)^2$.
The images of distinct degree-one inverse branches are disjoint.
Thus each marked degree-one component accounts for a distinct stationary root.
If several collisions occur at the same parameter, their multiplicities add. The
discriminant is the product of the discriminant of the moving factor, this squared
resultant, and the nonzero discriminant of $h$. This proves the zero count.

Let $N(x)$ count the zeros of $g$ of modulus at most $x$.
By Jensen's identity the mean in \eqref{eq:ascent} equals
$\int_0^{h_0}N(x)\,dx/(x\log 2)$. The logarithmic mean formula for
a linear factor also covers zeros on the circle. Substituting $r=x/2$
in the count and integrating the indicator for each $i$ gives
\eqref{eq:ascent}. Terms with $a_i=0$ vanish.
\end{proof}

The overlaps with the left and right eigenvectors enter this estimate
as a product, so both must be large for the same indices. We first
choose a sign vector with large overlaps with the left eigenvectors.
Restricting attention to the indices for which this succeeds, we then
choose a second sign vector with large overlaps with the right eigenvectors.

\begin{lemma}[Simultaneous choice of the two directions]\label{lem:signs}
Given an eigenbasis, its inverse, and nonnegative weights $\gamma_i$, two
vectors $u,v$ with entries $\pm1/n$ can be found in $O(n^2)$
operations such that $a_i\ge\|P_i\|/(4n^2)$ on indices carrying
at least $1/25$ of the total weight. For any $1\le m\le n$ specified indices,
suppose their unit left and right eigenvectors are known up to
individual phases with Euclidean error at most $1/(100n^2)$.
Let $s$ be the least power of two with $s^2\ge n$. In $O(nm)$ operations
one can choose $\bar u,\bar v$ with entries $\pm1/(2s)$ such that,
on at least $m/25$ indices, each overlap with the corresponding true
unit direction has modulus at least $1/(16s)$.
\end{lemma}
\begin{proof}
For independent uniform signs $\sigma$ and a nonzero complex column
$x$, define $Z=|x^*\sigma|^2/\|x\|^2$. Independence and symmetry give
$\mathbb E Z=1$ and $\mathbb E Z^2\le3$.
The polynomial $F(z)=z-z^2/6$ has maximum $3/2$ and is less than
$1/4$ on $[0,1/4)$. We fix the signs successively so that the conditional
expectation of $\sum_i\gamma_i F(Z_i)$ does not decrease. The final sum is at least
$\sum_i\gamma_i/2$, so indices with $Z_i\ge1/4$ carry at least
one fifth of the weight. To compute the conditional expectations,
let $m_0$ be the partial signed sum of a linear form and let $c_j$
be its remaining coefficients. Its fourth moment is
\[
 |m_0|^4+4|m_0|^2S+2\Re(\bar m_0^2 T)+2S^2+|T|^2-2U,
 \qquad S=\sum|c_j|^2,\quad T=\sum c_j^2,\quad U=\sum|c_j|^4 .
\]
The second moment is $|m_0|^2+S$. Updating the partial sum and the
three sums over the remaining coefficients costs constant work for each
linear form and sign. We apply this construction first to the left
eigenvectors and then to the right eigenvectors at the retained indices.
Dividing each sign vector by $n$ and using $P_i=x_iw_i^*$ proves the
first assertion.

For the second assertion, we apply the same two sign choices to the
$m$ pairs of approximate eigenvectors with equal weights. At each retained
index, the unscaled signed sum for either approximate vector has modulus
at least half its norm. That norm is at least $1-1/(100n^2)$, and replacing
the approximation by the true vector changes the
modulus by at most $\sqrt n/(100n^2)$. Division by $2s$ gives the
claimed lower bound. Each sign decision involves at most $m$ linear forms.
\end{proof}

\section{Discriminants of pencils and approximation of scalar roots}
\label{sec:scalar-batch}

Each rank-one update requires a choice of parameter on a circle.
We make this choice from the discriminant polynomial, allowing for
the cross term in the squared Frobenius norm of the perturbation.
To recover the eigenvalues after an update, we also need the roots of
the resulting characteristic polynomial to a prescribed accuracy.

All operation counts use the real-field model of
Section~\ref{subsec:main-model}, including comparisons when needed.
In purely algebraic operation counts,
$\widetilde O(m)$ denotes $O(m\log^c(m+2))$ for an absolute $c$.
The additional dependence on accuracy is specified in each numerical
subroutine. We use $\operatorname{polylog}(x)$ for a fixed power of
$1+\log x$, for $x\ge1$. The notation $\|p\|_1$ denotes the
sum of the moduli of the coefficients of $p$. Quasi-linear polynomial
multiplication over the coefficient field is available \cite{ck}.
Product and remainder trees then give quasi-linear multipoint evaluation
and interpolation at distinct nodes.

Newton iteration doubles the order of accuracy of a truncated reciprocal at each step.
Integrating $F'/F$ gives $\log F$ for $F(0)=1$, and Newton iteration on
$\log F=g$ gives the exponential for $g(0)=0$
\cite[Section~5]{brent-kung}. Polynomial translation reduces to a
convolution after factorial scaling. The extended polynomial gcd can be
computed by the half-gcd algorithm in $O(M(m)\log(m+2))$ operations,
where $M(m)$ is the
multiplication cost \cite[Chapters 9--11]{mca}. Each construction uses
field operations alone.

\subsection{The discriminant along a rank-one line}

The characteristic polynomial of a rank-one perturbation is affine in
the parameter. This allows us to obtain the resultants in its
discriminant from the traces of a multiplication operator.

\begin{lemma}[A resultant with one scalar parameter]\label{lem:paper-resultant}
Let $W\ne0$ have degree at most $2n$, and let $f,v$ have degrees at most $n$, with at least one nonzero. Up to multiplication by a fixed nonzero scalar, the resultant $\operatorname{Res}_z(W,f-tv)$ in the indeterminate $t$ can be computed in $\widetilde O(n^{(1+\omega)/2})$ operations, where $O(n^\omega)$ operations suffice for square matrix multiplication and $2\le\omega\le3$.
\end{lemma}
\begin{proof}
If $W$ is constant, the resultant is a nonzero constant. For
nonconstant $W$, a gcd computation determines whether $W,f,v$ have a
common zero, in which case the resultant vanishes identically. If they
have no common zero and $v=0$, the resultant is a nonzero constant.
In the remaining case, we normalise $W$ to be monic and define
\[
 M=W/\gcd(W,v^{\deg W+1}\bmod W)
\]
using the monic gcd. The removed roots, with their multiplicities,
contribute a nonzero constant, so we may return one if $M=1$.
Otherwise, $v$ is invertible in $\mathcal R=\mathbb C[z]/(M)$,
and we define $s=fv^{-1}\bmod M$. At each root $\zeta$ of $M$,
the operator of multiplication by $s$ on the local algebra is the
scalar $f(\zeta)/v(\zeta)$ plus a nilpotent operator. Its characteristic
polynomial, taken over all local factors, is therefore proportional to
the resultant as a polynomial in $t$.

The required traces can be computed by grouping the powers as in the
composition algorithm of Brent and Kung~\cite[Section~2.1]{brent-kung}.
Let $e=\deg M$, and choose the least integer $k$ with $k^2>e$.
Successive multiplication computes $s^a$ and $s^{jk}$ modulo $M$ for
$0\le a,j<k$ in $\widetilde O(ek)$ operations. For each $j$, we then
compute the first $e$ Laurent coefficients at infinity of
\[
 \frac{(M's^{jk})\bmod M}{M}.
\]
The coefficient of $z^{-i-1}$ is $\operatorname{tr}_{\mathcal R}(z^i s^{jk})$. Indeed, $M'/M$ has simple poles at the distinct roots, with their multiplicities as residues. Multiplication by a representative of $s^{jk}$ scales the residues by its values, and subtracting the polynomial part preserves all these Laurent coefficients. Computing the coefficients takes $\widetilde O(e)$ operations for each $j$.
Multiplying the resulting $k$-by-$e$ matrix by the $e$-by-$k$ matrix
of coefficients of the powers $s^a$ gives all traces
$\operatorname{tr}(s^i)$ for $1\le i\le e$. Splitting the inner
dimension into blocks of length $k$ bounds this multiplication cost
by $O((1+e/k)k^\omega)$. The reversed monic characteristic polynomial is
\[
 \exp\left(-\sum_{i=1}^e \frac{\operatorname{tr}(s^i)}{i}u^i\right)\bmod u^{e+1}.
\]
This identity follows by expanding the logarithm of the product over eigenvalues with multiplicity.
\end{proof}

\begin{proposition}[Discriminant and parameter choice]\label{prop:paper-phase}
Let $p$ be monic of degree $n\ge2$ with simple roots, let $\deg q<n$, and define $g(t)=\operatorname{disc}(p-tq)$. For $h_0>0$ one can compute $g$ up to nonzero scale and find $|t_0|=h_0$ with $|g(t_0)|\ge0.98\max_{|t|=h_0}|g(t)|$. Given at most $n$ further nodes $\xi$ with $p(\xi)\ne0$, one can choose $t=t_0(1+j/(1000n^4))$, $0\le j\le2n$, for which
\[
 |g(t)|\ge0.9\max_{|t|=h_0}|g(t)|,\qquad
 |p(\xi)-tq(\xi)|\ge |p(\xi)|/(10000n^4)\quad\text{for every such }\xi.
\]
The cost is $\widetilde O(n^{(1+\omega)/2})$. Strassen's bound \cite{strassen} gives $\omega=\log_2 7$ and $(1+\omega)/2<1.904$.
\end{proposition}
\begin{proof}
If $q=0$, the discriminant is a nonzero constant and any $t_0$ on the circle is admissible. Otherwise, let $h=\gcd(p,q)$ be monic, and define $f=p/h$ and $v=q/h$. The product formula gives
\[
 g(t)=\operatorname{disc}(h)\operatorname{disc}(f-tv)\operatorname{Res}(h,f-tv)^2,
\]
where the constant or linear discriminant is one. The first factor is nonzero because $p$ is simple. For $W=f'v-fv'$ the strict degree inequality $\deg f>\deg v$ gives $W\ne0$. Since $W=(f-tv)'v-(f-tv)v'$, reduction modulo the monic $f-tv$ yields
\[
 \operatorname{Res}(f-tv,W)=(-1)^{\deg f(\deg f-1)/2}\operatorname{disc}(f-tv)\operatorname{Res}(f,v).
\]
The final resultant is a nonzero constant. Two applications of
Lemma~\ref{lem:paper-resultant} therefore compute the full parameter
polynomial up to a single nonzero scalar, including all intersections
with stationary roots. Its scale, when needed, is determined by evaluating
the ordinary discriminant at $t=0$.

Homogeneity of the discriminant in the coefficients gives
$D=\deg g\le2n-2$. The parametrisation $(1+is)/(1-is)$,
$0\le s\le1$, together with its rotations by multiples of $\pi/2$,
covers the unit circle with speed at most two. We sample these arcs at
rational steps $1/(100(D+1))$ and multiply each sample by $h_0$.
To bound the loss between samples, consider a polynomial $P$ of degree
at most $D\ge1$ and let
$M_P=\max_{|z|=1}|P(z)|$. If $|a|>M_P$, Rouch\'e's theorem places
all $D$ roots of $az^D-P(z)$ in the open unit disc. By Gauss--Lucas its
derivative is nonzero on the unit circle. Taking
$a=P'(z)/(Dz^{D-1})$ at a boundary point therefore gives
$|P'(z)|\le DM_P$. Scaling the circle gives Bernstein's derivative
inequality at any positive radius. Consequently, one of our samples
attains a modulus within $0.02M_0$ of the maximum, where
$M_0=\max_{|t|=h_0}|g(t)|$. Multipoint evaluation and comparisons
of squared moduli select such a sample $t_0$. For $D=0$, we take $t_0=h_0$.

Set $u=2n/(1000n^4)$. The exterior polynomial bound at radius $h_0(1+u)$ is $(1+u)^D M_0$, by the maximum principle applied to the reversed polynomial. Bernstein's inequality there bounds the derivative on the expanded disc by $D(1+u)^DM_0/h_0$. For each trial $j$, the decrease in modulus is at most $uD(1+u)^DM_0<0.08M_0$. For each prescribed node, let $a=p(\xi)$, $b=t_0q(\xi)$. If $|b|\le|a|/2$, every trial satisfies the required bound. Otherwise the successive values of $a-b(1+j/(1000n^4))$ are more than $|a|/(2000n^4)$ apart. At most one trial fails the bound. It suffices to test a closest integer in $[0,2n]$ to $1000n^4\Re(a/b-1)$, found by comparisons. At most $n$ indices are excluded. All work following the construction of $g$ costs $\widetilde O(n)$.
\end{proof}

The squared Frobenius norm of the accumulated perturbation contributes
a linear term to the logarithmic objective on the circle. Exponentiating
gives the weight in the next proposition.

\begin{proposition}[Phase selection with an exponential weight]\label{prop:weighted-phase}
Let $p(x)\in\mathbb C[x]$ be monic of degree $n\ge2$ with simple
roots, let $q(x)\in\mathbb C[x]$ have degree less than $n$, and let
$r>0$. Define $g(z)=\operatorname{disc}_x(p(x)-rzq(x))$.
Let $c\in\mathbb C$ and suppose we are given $c_0$ satisfying
$|c_0-c|\le1/100$. Let $A$ be an integer power of two in $[1+|\Re c_0|+|\Im c_0|,\,2(1+|\Re c_0|+|\Im c_0|)]$,
and $J=2n+64A$. Given at most $n$ nodes $\xi$ with $p(\xi)\ne0$,
one can compute $z$ such that
\[
 1\le |z|\le 1+\frac{1}{500J},\qquad
 |g(z)e^{-cz}|\ge 0.9\max_{|w|=1}|g(w)e^{-cw}|,\qquad
 |p(\xi)-rzq(\xi)|\ge\frac{|p(\xi)|}{2^{16}nJ}
\]
at every given node. The cost is
$\widetilde O(n^{(1+\omega)/2}+J)$.
\end{proposition}
\begin{proof}
Form the degree-$64A$ Taylor polynomial $T$ of $e^{-c_0z}$ by the
coefficient recurrence. For $|z|\le2$ and $l=64A+1$, its relative
error is at most $e^{4A}(2eA/l)^l<1/100$, by the exponential
series remainder and $l!\ge(l/e)^l$. Thus the modulus of $gT$ is
between $0.99e^{-0.02}$ and $1.01e^{0.02}$ times that of
$ge^{-cz}$ throughout this disc. Proposition~\ref{prop:paper-phase}
computes $g$ up to scale. The product $gT$ has degree at most $J$.
The rational circle samples and multipoint evaluation in the proof
of that proposition give a point $z_0$ attaining at least $0.98$ of
its maximum modulus $M$ on the circle.

We now test the points $z_j=z_0(1+j/(1000nJ))$, $0\le j\le2n$.
On the expanded disc of radius $1+1/(500J)$ the polynomial has
modulus at most $(1+1/(500J))^J M$. Bernstein's inequality on its
boundary and the maximum principle for the derivative give
$|(gT)(z_j)-(gT)(z_0)|<0.003M$. Since
$0.977\,(0.99/1.01)e^{-0.04}>0.9$, every trial has the required
weighted modulus. For a given node, let $b=rz_0q(\xi)$. If
$|b|\le |p(\xi)|/2$ all trials satisfy the required lower bound. Otherwise consecutive values
of $p(\xi)-b(1+j/(1000nJ))$ differ by more than
$|p(\xi)|/(2000nJ)$, so at most one trial fails. A closest integer
in $\{0,\ldots,2n\}$ to $1000nJ\Re(p(\xi)/b-1)$ can be found by
binary search when $b\ne0$. We test the corresponding trial and
exclude it only if it fails the lower bound. Since each node excludes
at most one index, at least one trial remains. Polynomial multiplication, sampling and the
node tests cost $\widetilde O(J)$ after $g$ is known.
\end{proof}

\subsection{Root approximation in the real-field model}

Root-radius estimates and recursive splitting circles are central to
the algorithms of Neff and Reif~\cite{neff-reif} and
Pan~\cite[Sections~2.2--2.5]{pan}. We need a real-field operation count with explicit dependence on the
requested accuracy. We control the error in coefficient norm and use
the following lemma to match the roots, counting multiplicities.

\begin{lemma}[Matching roots from a coefficient residual]\label{lem:paper-match}
If $p,q$ are monic of degree $m\ge1$, the roots of $p$ lie in the unit disc, and $\|q-p\|_1\le e<2^{-m}$, the roots of $q$ can be matched with those of $p$ to distance at most $4m e^{1/m}$, including multiplicities.
\end{lemma}
\begin{proof}
On $|z|=2$ we have $|p(z)|\ge1$ and $|q(z)-p(z)|\le2^me<1$. Rouch\'e's theorem keeps the zeros of $p_t=p+t(q-p)$ inside this circle for $0\le t\le1$. At any such zero, $|p(z)|\le2^me$, so its distance to the nearest root of $p$ is at most $2e^{1/m}$. We enlarge this radius by an arbitrary $\varepsilon>0$. Each connected component of the resulting union of discs has diameter at most $2m(2e^{1/m}+\varepsilon)$. The roots move continuously as a multiset, so the count in each component stays constant. Matching within each component and letting $\varepsilon$ decrease to zero along a subsequence with fixed matching proves the assertion.
\end{proof}

\begin{theorem}[Scalar root approximation]\label{thm:paper-roots}
For a monic polynomial $p$ of degree $m\ge1$ with roots in the closed unit disc and an integer $b\ge2$, there is a deterministic algorithm producing $m$ field elements $z_i$ such that
\[
 \left\|p-\prod_{i=1}^m(z-z_i)\right\|_1\le2^{-b}.
\]
It costs $O(m^{1.9}\operatorname{polylog}(m(b+2)))$ real-field operations and comparisons. In particular, $b\ge m(\mu+\lceil\log_2(4m)\rceil)$ gives absolute matched root accuracy $2^{-\mu}$ for $\mu\ge0$.
\end{theorem}

Our algorithm factors the polynomial recursively. Enclosures of the
root distances guide the search for a circle that separates the roots
into two groups, or for a small cluster. Either outcome gives two
approximate factors. At the start of each recursive call, we increase
$b$ to at least $2m\lceil\log_2(m+2)\rceil+10$ and define
$B=b+20m+20$, $\eta=2^{-B}/8$. Integer exponents, dyadic powers, and the
precision bounds below are obtained by doubling, squaring, and binary
search on integers.

The next four lemmas concern a monic polynomial $p$ of degree
$m\ge2$ whose roots lie in the closed unit disc, with $B$ and $\eta$
defined above.

\begin{lemma}[Computing enclosures of root radii]\label{lem:paper-radii}
Let $C$ be a field element with $|C|\le16$, and define $\rho=2^{-35}/m^2$. Denote the roots of $p$, ordered by distance from $C$ and counted
with multiplicity, by $\zeta_1,\ldots,\zeta_m$.
In $m\operatorname{polylog}(m(B+2))$ operations one can enclose the ordered values $|\zeta_i-C|$ for the true roots in intervals $[L_i,U_i]$. An initial segment has enclosures $[0,(1+\rho)\eta]$. Every remaining enclosure has $L_i>0$ and $U_i/L_i\le1+4\rho$. Every root at distance less than $\eta/(1+\rho)$ from $C$ has an enclosure beginning at zero. The two endpoint lists may both be taken in nondecreasing order.
\end{lemma}
\begin{proof}
Let $N$ be the least power of two at least $2^{38}m^3$.
We translate $p$ by $C$ and apply Graeffe squaring, using its even and
odd parts, until the roots are the $N$th powers of the translated roots.
After removing zero roots by exact division, we use the nonzero
coefficients in the hull construction below. If the polynomial is a
monomial, all radii are zero. Otherwise, let $p_1$ be the remaining
polynomial, with coefficients $a_j$ and nonzero roots $r_i$ ordered by
nondecreasing modulus. Consider
\[
 H(t)=\max_{a_j\ne0}(\log|a_j|+jt),\qquad
 J(t)=\log|a_{\mathrm{lead}}|+\sum_i\max(t,\log|r_i|).
\]
The circle mean of $\log|p_1(e^te^{i\theta})|$ equals $J(t)$ by the logarithmic mean formula for each linear factor. Hence $J\le H+\log(m+1)$. The elementary symmetric coefficient formulas give $H\le J+m\log2$. Between the $i$th breakpoints of these functions, counting each slope
jump with its multiplicity, the function whose breakpoint occurs earlier
has slope at least $i$, while the other has slope at most $i-1$.
Their breakpoints are therefore separated by at most
$m\log2+\log(m+1)$ in $t$. Since the $i$th breakpoint of $J$ is
$\log|r_i|$, taking $N$th roots gives estimates of all original positive
radii within a multiplicative factor $1+\rho$.
Indeed, the logarithmic error is at most
$(m\log 2+\log(m+1))/N\le2m\log 2/N$, whose exponential is less than $1+\rho$.

We locate the breakpoints using comparisons of integer powers.
A monotone scan over the already ordered abscissae constructs the upper
hull of $(j,\log|a_j|^2)$ for the nonzero coefficients. Each orientation
test compares integer powers, computed by repeated squaring. There
are $O(m)$ tests, costing $O(m\log(m+2))$ operations in total. An edge from $i$ to $i+k$ gives the ordered estimate
$(|a_i|^2/|a_{i+k}|^2)^{1/(2kN)}$ with multiplicity $k$.
Comparison of integer powers determines whether this estimate is at
most $\eta$. In that case, or when the true radius is zero, we use the
displayed enclosure beginning at zero. Otherwise, binary search on the
exponents finds an interval between consecutive powers of two that
contains the estimate and lies between $\eta$ and the fixed upper bound
$64$. We bisect this interval until its relative width is at most
$\rho/4$, then divide its lower endpoint by $1+\rho$ and multiply its
upper endpoint by $1+\rho$. The resulting width ratio is at most $1+4\rho$. One may use fixed subdivisions of relative mesh at most $\rho/4$ in each such interval, with monotone bracketing at shared endpoints. For a positive lower endpoint, the upper bracket for the estimate is at least $\eta$, so the enlarged upper endpoint is at least $(1+\rho)\eta$. Thus both endpoint lists are monotone. The cost of this search is polynomial in $\log(m(B+2))$, including the powers used in comparisons.
\end{proof}

Roots of higher derivatives provide centres for splitting circles in
Neff and Reif~\cite{neff-reif} and Pan~\cite[Section~2.4]{pan}.
For each disc containing enough roots, an elementary estimate gives a
nearby derivative root. Applying the estimate to a disc of least radius
allows us to choose the same derivative root for all such discs.

\begin{lemma}[A common derivative centre]\label{lem:paper-derivative}
Let $p$ have degree $m$, let $l$ be an integer with $m/2\le l<m$, and define $d=m-l$. Some zero of $p^{(l)}$ lies in $D(Y,(2d+2)R)$ for every closed disc $D(Y,R)$ containing at least $l+1$ zeros of $p$, counted with multiplicity.
\end{lemma}
\begin{proof}
First fix one such disc of positive radius. Translation and scaling
reduce it to the unit disc. Choose $l+1$ of its roots
$\alpha_0,\ldots,\alpha_l$, counting multiplicities, and express
$p(z)=\sum a_jz^j$. The quotient by $\prod_{i=0}^l(z-\alpha_i)$ is a polynomial, so its Laurent coefficient at $z^{-1}$ vanishes. Consequently
\[
 0=a_l+\sum_{r=1}^{d}a_{l+r}h_r,\qquad
 h_r=[u^r]\prod_{i=0}^l(1-\alpha_i u)^{-1},\qquad
 |h_r|\le \binom{l+r}{r}.
\]
If every root of $p^{(l)}$ had modulus greater than $2d$, then $a_l\ne0$. Factoring $p^{(l)}(z)/p^{(l)}(0)$ into $d$ factors with linear coefficients of modulus less than $1/(2d)$ would give
\[
 |a_{l+r}|\binom{l+r}{r}\le |a_l|\binom{d}{r}(2d)^{-r}.
\]
The sum over $r\ge1$ would then have modulus at most $|a_l|((1+1/(2d))^d-1)<|a_l|$, a contradiction. For radius zero the same local assertion follows from the multiplicity at the centre.

We now choose a disc of minimum radius $R_*$ containing at least
$l+1$ roots. Such a disc exists by finiteness of the root multiset and
compactness. The preceding argument gives a derivative root in the disc
with the same centre and radius $2dR_*$. We fix this derivative root.
Every other disc in the statement shares a root with the minimum disc
because $2(l+1)>m$. The distance between their centres is at most $R+R_*$, and $R_*\le R$, so this derivative root lies within $(2d+2)R$ of every such centre.
\end{proof}

\begin{lemma}[Finding a separating annulus]\label{lem:paper-annulus}
For $m>200$, define $g=\lceil m/50\rceil$, $h=m-g$, $l=\lfloor17m/20\rfloor$, and $s=m+3$. Given approximants to the roots of $p^{(l)}$ matched to error $a=\eta/(100s)$, one can find either of the following in $m\operatorname{polylog}(m(B+2))$ operations.
\begin{enumerate}
\item A centre $C$ and circle of radius $r$ with exactly $k\in[g,h]$ roots inside it, such that the inner roots have distance at most $(1-\sigma)r$ from $C$ and the outer roots have distance at least $(1+\sigma)r$, where $\sigma\ge c/m$ for an absolute rational $c>0$.
\item A centre $C$ with a known number $k>h$ of roots in a disc of radius less than $4\eta$.
\end{enumerate}
In both cases $|C|,r\le16$, where in the second case $r$ denotes the disc radius.
\end{lemma}
\begin{proof}
Let $\psi=1+2^{-20}/m$. At each candidate centre, we apply
Lemma~\ref{lem:paper-radii}. If at least $g$ enclosures begin at zero,
we consider the $m+1$ closed annuli with inner and outer radii
$2\eta\psi^{3j}$ and $2\eta\psi^{3j+1}$, $0\le j\le m$.
The enclosures beginning at zero lie inside the inner boundary of every
annulus. Each remaining enclosure meets at most one annulus, since the
ratio separating successive annuli is $\psi^2>1+4\rho$. Thus at least
one annulus meets no enclosure. Sorting and scanning the endpoints, or using binary
search to identify the annulus met by each enclosure, finds such an
annulus and counts its interior roots. If this count is at most $h$,
the circle halfway between the two boundaries gives the first outcome.
Otherwise, the inner radius is below $4\eta$, giving the second outcome.
If fewer than $g$ enclosures begin at zero, any gap
$L_{i+1}\ge\psi U_i$ with $g\le i\le h$ gives the first outcome
using the circle of radius $(U_i+L_{i+1})/2$. Both midpoint constructions
give a separation ratio of at least $(\psi-1)/(\psi+1)$.

In the remaining case the shell $[L_g,U_{h+1}]$ contains at least $m-2g+2$ roots and its ratio of outer to inner radii is at most
\[
 (\psi(1+4\rho))^{m+1}\le1+4\cdot2^{-20}.
\]
Let $r_0$ be its outer radius and $C$ its centre. We repeat the
radius-enclosure computation at $C+2r_0$ and $C+2ir_0$. If either
computation gives one of the stated outcomes, the construction is
complete. Otherwise, the three shells share at least
$m-6g+6\ge22m/25>l$ roots. Let their squared radius intervals be
$[R_j-\Delta_j,R_j+\Delta_j]=[a_j^2,b_j^2]$ for $j=0,1,2$,
respectively. Subtracting the squared-distance inequalities in the
real and imaginary directions places their common intersection in the
disc with centre and radius
\[
 Y=C+\frac{4r_0^2+R_0-R_1}{4r_0}
       +i\frac{4r_0^2+R_0-R_2}{4r_0},\qquad
 R=\frac{2\Delta_0+\Delta_1+\Delta_2}{4r_0}.
\]
The two coordinate errors are bounded by
$(\Delta_0+\Delta_1)/(4r_0)$ and $(\Delta_0+\Delta_2)/(4r_0)$.
Their sum bounds the Euclidean error and gives the displayed radius.
A common point is at most $3r_0$ from either offset centre, so $b_j\le4r_0$ by the relative widths. These widths give $R\le r_0/100$. Since a common root is between $r_0/(1+4\cdot2^{-20})$ and $r_0$ from $C$, we also have $r_0/2\le|Y-C|\le2r_0$.

We continue with $Y$ as the next centre. If the three new shells
again determine a disc, the first shell contains the $g$th nearest root
distance, which is at most $R$. Its relative width, together with the
possible enclosure error $1+4\rho$ for that distance, bounds its outer
radius by $r_0/50$. If $r_0\le50\eta$, then $R\le\eta/2$ and at
least $g$ enclosures at $Y$ begin at zero. The next computation therefore
gives one of the two stated outcomes.

For the starting centre, we use the coordinatewise median $C_0$ of
the retained approximants to the derivative roots, counting
multiplicities. Since $2(m-l)+2\le s$, Lemma~\ref{lem:paper-derivative}
gives a single derivative root that lies in each of the discs determined
by three shells after its radius is multiplied by $s$. We retain an
approximant matched to this root throughout the construction.
Starting at $C_0$, we generate up to $J=2+\lceil\log_2(100s)\rceil$
consecutive discs determined by three shells, stopping as soon as either
outcome of the lemma occurs. If all $J$ discs are produced, let $Y,R$ be the final centre and radius and let $r_0$ now denote the outer radius of the first shell. Then $r_0>50\eta$. The first displacement has size at least $r_0/2$ and subsequent displacements have sum at most $2r_0/49$, so $|Y-C_0|>r_0/3$. Geometric decay gives $sR\le r_0/100$. The disc $D(Y,sR+a)$ contains the prescribed matched approximant and has radius below $r_0/50$. It lies strictly to one side of either the real or the imaginary median coordinate. Removing every approximant outside this disc leaves a nonempty set of at most half the original size, even with tied coordinates. After $O(\log(m+2))$ such chains, one of the two outcomes must occur.

It remains to bound the centres at which we apply Lemma~\ref{lem:paper-radii}. All derivative roots have modulus at most one by Gauss--Lucas, and their matched approximants have modulus less than two. For a coordinatewise median, $|C_0|\le3$. The outer radius of the initial shell at this centre is at most $(1+4\rho)(|C_0|+1)$, because its upper endpoint bounds a radius whose enclosure has positive lower endpoint. Thus the two offset centres have modulus below $12$. A subsequent centre obtained from three shells lies within $r_0/100$ of a unit-disc root. Along this sequence the shell radii decrease, so the same bounds apply throughout. The radius of any midpoint circle is bounded by $(1+4\rho)(1+|C|)$ or by $4\eta$, and is therefore below $16$. All comparisons use squared moduli. The number of radius-enclosure computations is $O(\log^2(m+2))$.
\end{proof}

A separating annulus allows us to approximate the power sums of the inner
roots. For a small cluster we can use its centre for each inner root.

\begin{lemma}[Recovering the two factors]\label{lem:paper-recover}
A circle from the first outcome of Lemma~\ref{lem:paper-annulus} gives monic factors $f_o,g_o$ of degrees $k,m-k$ with $\|p-f_og_o\|_1<2^{-b-3}$. In the second outcome one may take $f_o=(z-C)^k$. Every root of the computed factors is in $|z|<2$. The cost is $m\operatorname{polylog}(m(B+2))$.
\end{lemma}
\begin{proof}
For a separated circle, transform $p$ to the monic polynomial $S(y)=r^{-m}p(C+ry)$. Let $y_i$ denote its roots. For a power of two $N\ge2m$ compute
\[
 h=(1+y^N)^{-1}\bmod S,\qquad
 F_j=[y^{-j-1}]\frac{(S'h)\bmod S}{S},\quad 1\le j\le k.
\]
The separating annulus guarantees invertibility. Repeated modular squaring, extended gcd and reversed division cost $\widetilde O(m\log N)$. The same partial-fraction computation as in Lemma~\ref{lem:paper-resultant} gives $F_j=\sum_i y_i^j/(1+y_i^N)$ including multiplicities. For $|y_i|\le1-\sigma$, the absolute difference from $y_i^j$ is at most $2(1-\sigma)^N$. For $|y_i|\ge1+\sigma$ the contribution is at most $2(1+\sigma)^{-(N-m)}$, once $(1+\sigma)^{-N},(1-\sigma)^N\le1/2$. Take
\[
 b'=b+40m+40+\lceil\log_2(m+1)\rceil
\]
and choose $N$ as the next power of two above a sufficiently large universal multiple of $m(m+b'+\log(m+2))/c$, with $c$ the separation constant. All $k$ sums then approximate the true inner power sums to absolute error $e=2^{-b'}$.

We recover the approximate monic inner factor from its reversed polynomial
\[
 \exp\!\left(-\sum_{j=1}^k F_j u^j/j\right)\bmod u^{k+1}.
\]
The exact reversed inner factor has coefficient norm at most $2^m$. The perturbation to the exponent has norm at most $me$, so the factor error is at most $2^m(\exp(me)-1)\le2^{m+1}me$. After reversing the list of $k+1$ coefficients, we return to the original
variable by the map $f\mapsto r^k f((z-C)/r)$. Since $r,|C|\le16$,
this map amplifies the coefficient norm by at most $33^m$. The resulting factor $f_o$ in the original variable is monic and differs from the true monic factor $F$ by at most $e_0=2^{-b-20m-20}$.

Let $G$ denote the true monic outer factor and take $g_o$ to be the quotient in Euclidean division of $p$ by $f_o$. We have $\|G\|_1\le2^m$. The inverse of the reversed $F$ modulo $u^{m+1}$ has norm at most $\binom{2m}{m}\le4^m$, by expanding the reciprocal linear factors at unit-disc roots. By a Neumann series, the corresponding inverse of the reversed $f_o$ has norm at most $2\cdot4^m$. Reversed division in $p=f_oG+(F-f_o)G$ thus gives $\|g_o-G\|_1\le2\,8^m e_0$. The triangle inequality and $\|f_o\|_1\le2^m+e_0$ now give $\|p-f_og_o\|_1<2^{-b-3}$.

In the second case, we set $w=z-C$ and divide $p(C+w)$ by $w^k$.
The quotient $Q_w$ consists of the coefficients of degrees $k$ and above,
with their degrees decreased by $k$. The translated true inner factor $F$ satisfies $\|F-w^k\|_1\le(1+r)^k-1\le2mr$, where $r<4\eta<2^{-B}$. The centre lies within $r$ of a unit-disc root, hence $|C|<2$. The translated outer factor $G$ therefore has norm at most $4^m$, giving $\|Q_w-G\|_1\le2mr4^m$. The residual before translating back is at most $4mr4^m$. Taking $g_o(z)=Q_w(z-C)$ gives residual at most $4mr\,12^m<2^{-b-3}$. Finally, in both cases Rouch\'e's estimate in Lemma~\ref{lem:paper-match} places all zeros of $f_og_o$ inside the circle of radius two.
\end{proof}

\begin{proof}[Proof of Theorem~\ref{thm:paper-roots}]
For $m>200$, set $l=\lfloor17m/20\rfloor$, $d=m-l$, and $s=m+3$
as in Lemma~\ref{lem:paper-annulus}. We first recursively approximate
the roots of the monic normalisation of $p^{(l)}$, whose degree is
$d\le0.155m$. They are in the unit disc by Gauss--Lucas.
Lemma~\ref{lem:paper-match} provides the required matched accuracy
$a=\eta/(100s)$ when this recursive approximation uses the precision
exponent $d(\lceil\log_2(4d/a)\rceil+1)$ for the coefficient residual.
Lemmas~\ref{lem:paper-annulus} and \ref{lem:paper-recover} then give two
approximate factors.

For $2\le m\le200$, we use $g=1$, $h=m-1$ and the arithmetic mean
of the roots as the first centre. This mean is obtained from the
subleading coefficient. If the same three radius-enclosure computations
gave neither a separating circle nor a cluster in the smaller disc,
each shell would contain all roots. Their common disc $D(Y,R)$ would
therefore contain all roots and their mean $C$, contradicting
$|Y-C|\ge r_0/2$ and $R\le r_0/100$. Hence the derivative roots
need not be approximated for these degrees. A polynomial of degree one
gives its root by a field operation, and an empty factor gives the empty list.

Before recursively factoring a remaining factor of degree $j$, we
substitute $z=2u$ and divide by $2^j$ so that its roots lie in the
unit disc. We use a precision exponent of at least $b+20m+20$ for the
coefficient residual, together with the lower bound imposed at the
start of each call. For a cluster in
the smaller disc, we return $k$ copies of $C$ and recurse only on
the other factor. Reversing this change of variable amplifies the coefficient norm by a factor of at most $2^m$. Each original factor has coefficient norm at most $3^m$ since its roots lie inside $|z|<2$. The product of the returned linear factors therefore introduces additional error at most $2\cdot3^m\varepsilon+\varepsilon^2$, where $\varepsilon=2^m2^{-b-20m-20}$. Together with the initial residual this gives the asserted $2^{-b}$.

The precision exponent required in each recursive call is at most $C_2m(b+m+\log(m+2))$ for an absolute $C_2$. For degrees above $200$, each polynomial in a recursive call has degree at most $0.98m$ and the sum of the $1.9$th powers of the degrees in those calls is at most
\[
 m^{1.9}(0.155^{1.9}+0.98^{1.9}+0.02^{1.9})<0.993m^{1.9}.
\]
Indeed, the fraction of roots in either factor ranges from $0.02$ to $0.98$, and convexity maximises the sum at either endpoint. In the second outcome, recursion is needed only for the remaining
factor of degree $<0.02m$ and for the derivative polynomial. The logarithm of the precision exponent on any path grows by $O(\log^2(m+2))$ over its initial logarithm before reaching bounded degree. The recursion for the finitely many remaining degrees adds only a fixed constant to the depth and a fixed multiplier to that precision exponent. Thus the local polylogarithmic factors are uniformly bounded by a fixed polynomial in $\log(m(b+2))$. Summing the costs over the recursion tree as a geometric series proves
the operation bound. The matched root assertion follows from Lemma~\ref{lem:paper-match}.
\end{proof}

\section{Initialisation with a separated spectrum}\label{sec:seed}

The discriminant iteration starts from a matrix with a separated spectrum
and a known eigensystem. A rank-one correction makes the approximate
roots of a characteristic polynomial exact eigenvalues. We then increase
the gaps by shifting the diagonal entries in triangular form.

\begin{proposition}[Initial separated eigensystem]\label{prop:seed}
For $n\ge2$, $\|A\|\le1$, and $0<\delta<1/2$, one can compute an exact
diagonalisation of a matrix $B_0$ satisfying
\[
  \|B_0-A\|<\delta/4,\qquad
  \operatorname{gap}(B_0)\ge\delta/(8n),
\]
in $O(n^3\operatorname{polylog}(n/\delta))$ real-field operations and comparisons.
\end{proposition}
\begin{proof}
We first reduce $A$ to upper Hessenberg form by conjugations on
adjacent coordinates. The second entry of a nonzero pair $(a,b)$ is
eliminated by multiplying on the left by
\[
 \frac1s\begin{pmatrix}\bar a&\bar b\\-b&a\end{pmatrix},
 \qquad
 \sqrt{|a|^2+|b|^2}\le s\le
 (1+1/(100n^4))\sqrt{|a|^2+|b|^2},
\]
and on the right by its exact inverse. We compute $s$ by the component
scaling and bisection of Section~\ref{sec:foundations}, omitting zero
pairs. There are fewer than $n^2$ conjugations, each with condition
number at most $1+1/(100n^4)$, and their product has condition number
below $1.02$.

Let $\eta=\delta/128$. Replacing each subdiagonal entry of modulus
less than $\eta$ by $\eta$ gives a correction supported on a single
subdiagonal, with norm at most $2\eta$. The resulting matrix $H$
satisfies $\|H\|<2$, and $e_1$ is cyclic for $H$.

Indeed, $v_j=H^je_1$ vanishes below coordinate $j+1$ and its
$(j+1)$st entry is the product of the first $j$ subdiagonal entries
for $0\le j<n$. The matrix $V=[v_0,\ldots,v_{n-1}]$ therefore
has determinant of modulus at least $\eta^{n(n-1)/2}$.
Computing $v_0,\ldots,v_n$ and expressing $v_n$ in this basis gives
the recurrence corresponding to the monic characteristic polynomial $p$.
The polynomial part of $p(z)\sum_{j\ge0}v_jz^{-j-1}$ is the column
$r(z)=\operatorname{adj}(zI-H)e_1$.
The coefficient matrix of $r$, in descending order, is $V$ times
a unit upper triangular Toeplitz matrix. Its determinant satisfies the
same lower bound. Its norm is at most $n^2 6^n$ by
$\|v_j\|\le2^j$ and the coefficient bound from the roots of $p$.
The cofactor formula therefore bounds the logarithm of its inverse norm by
$O(n^2\log(2n/\delta))$.

For a target accuracy $0<e\le1$, we apply
Theorem~\ref{thm:paper-roots} to $2^{-n}p(2z)$ with matched accuracy
$e/4$. This gives accuracy $e/2$ in the original variable. We order
the resulting nodes by real part and add $je/(2n)$ to the $j$th node,
$0\le j<n$. The resulting nodes $\mu_j$ are distinct and still match
the roots of $p$ to within $e$. Let $f(z)=\prod_j(z-\mu_j)$.
For $e\le1$, telescoping the products of linear factors bounds
$\|p-f\|_1\le n5^ne$. Solving the nonsingular coefficient system
$\ell^*r=p-f$ gives a rank-one correction
$H_1=H+e_1\ell^*$ whose characteristic polynomial is $f$, by the
determinant lemma. The preceding inverse bound shows that
$\|\ell\|\le\delta/32$ once $\log_2(1/e)$ is a sufficiently
large constant times $n^2\log(2n/\delta)$. One may double this
exponent until the squared norm meets the bound. All $r(\mu_j)$
are nonzero by invertibility of their coefficient matrix, and
$(zI-H_1)r(z)=f(z)e_1$ gives exact eigenvectors.

We apply Gram--Schmidt to these eigenvectors in increasing order of
the real parts of their eigenvalues. Exact projection coefficients and
approximate length scaling give an orthogonal basis of the invariant
flag, with condition number at most $1+1/(100n^4)$. In this basis,
$H_1$ is triangular. Increasing the $j$th diagonal entry by
$j\delta/(8n)$, $0\le j<n$, separates the real parts of successive
eigenvalues by at least $\delta/(8n)$. Triangular substitution then
gives all eigenvectors. The product of the condition numbers of the
Hessenberg and orthogonal coordinate changes is below $1.03$, so the
perturbation in the original coordinates is bounded by
\[
 1.03(2\eta+\delta/32+\delta/8)<\delta/4 .
\]
The Hessenberg reduction, Krylov sequence, coefficient solves, and
basis constructions take cubic work up to logarithmic factors.
The cost of root approximation depends only polylogarithmically on
the precision exponent, and the number of exponent doublings is logarithmic.
\end{proof}

\section{Resolvents of low-rank perturbations}
\label{sec:transfer-rebase}

We use the same reference eigenbasis for several successive rank-one
perturbations. In these coordinates, the determinant lemma reduces the
characteristic polynomial to a small determinant, and the resolvent gives
approximations to the left and right eigenvectors. A further small
correction makes the computed eigenvalue approximations exact.

\subsection{Reference coordinates and characteristic pencils}

Let $n\ge2$ and $B_0=X_0\Lambda W_0$, where $W_0=X_0^{-1}$ and
$\Lambda=\operatorname{diag}(\lambda_1,\dots,\lambda_n)$ has distinct
entries. The current matrix has the form $B=B_0+U K V^*$, where the
perturbation factors have at most $k$ columns. Their coordinates
$W_0U$ and $V^*X_0$ in the reference basis are available. Throughout this section, we assume that
$\|B_0\|,\|B\|\le2$. The notation $\widetilde O$ suppresses fixed
powers of $1+\log n$ in the algebraic subroutines.

The cost estimates use the deterministic multiplication theorem of
Bostan, Jeannerod, Mouilleron and Schost
\cite[Theorem 1.1]{structured}. For inversion we use positive definite
Schur complements. Suppose that the diagonal matrices $D_a,D_b$ have orders at most
$n$, distinct diagonal entries within each matrix, and disjoint spectra. If
$D_a C-C D_b=UV^*$ is given by generators of width $w\ge1$, the
product of $C$ with $w$ dense columns costs
$\widetilde O(nw^{\omega-1})$ field operations for $w\le n$.
We apply the theorem to the linear factors associated with the diagonal
entries. These are monic and pairwise coprime within each family, and
disjointness of the spectra makes the Sylvester operator invertible. The coefficient field has characteristic
zero, and the companion matrices of the linear factors are the scalar
nodes themselves. Transposition gives
the same bound for multiplication on the left. Padding or blocking
$b$ right-hand sides therefore gives the cost
$\widetilde O(n\max(b,w)w^{\omega-2})$. For a block of dimension at
most $w$, the same bound follows from ordinary dense multiplication
with square blocking.

\begin{lemma}[Characteristic pencils]\label{lem:batch-pencil}
Given the reference coordinates of the perturbation factors,
including $W_0u$ and $v^*X_0$ for a proposed update $suv^*$, the
characteristic polynomial of $B$ and the pencil for this update can
be formed in
$\widetilde O(n(k+1)^3)$ field operations.
\end{lemma}
\begin{proof}
Let $p_0(z)=\prod_i(z-\lambda_i)$. By the determinant lemma,
\[
 p_K(z)=p_0(z)\det\!\left(I-K V^*X_0
             (zI-\Lambda)^{-1}W_0 U\right).
\]
We evaluate this identity at $z=4+j/(n+1)$, $0\le j\le n$.
The $k^2$ entries of the small matrix are scalar Cauchy sums. Their
numerators over the common denominator $p_0$ are computed by a product
tree and evaluated at the external nodes by a remainder tree. This costs
$\widetilde O(nk^2)$. Exact pivoting computes the $n+1$ determinants
in $O(nk^3)$ operations, and interpolation recovers $p_K$. Finally,
the characteristic polynomial is affine in an additional rank-one
update $suv^*$. If $p$ and $p_1$ are the current polynomial and the polynomial at $s=1$,
respectively, then $q=p-p_1$ gives the pencil $p-sq$.
\end{proof}

The phase search also requires the discriminant of this pencil,
including collisions with stationary roots. Let $h_0=\gcd(p,q)$ and
factor $p=h_0p_1$, $q=h_0q_1$. Since the current spectrum is simple,
$h_0$ is square-free and coprime to $p_1$. The discriminant of $p-sq$
is
\[
 \operatorname{disc}(h_0)\operatorname{disc}(p_1-sq_1)
           \operatorname{Res}(h_0,p_1-sq_1)^2 .
\]
If $q\ne0$, Section~\ref{sec:scalar-batch} computes the second factor.
In the resultant, the roots of $\gcd(h_0,q_1)$ contribute a nonzero
constant. On the quotient by $h_0/\gcd(h_0,q_1)$, the element
$p_1/q_1$ is well defined. The power-projection calculation of
Section~\ref{sec:scalar-batch} gives its characteristic polynomial,
which is proportional to $\operatorname{Res}(h_0,p_1-sq_1)$.
Its square gives the last factor up to a nonzero constant. Evaluation at $s=0$
fixes the constant if it is required. If $q=0$, the full discriminant
is constant.

\subsection{Separation from the reference spectrum}

The structured matrix products require the current and reference spectra
to be disjoint. We first ensure this by a rank-one perturbation.
A small adjustment to the parameter radius preserves disjointness at
subsequent steps.

\begin{lemma}[Separation from the reference spectrum]\label{lem:old-poles}
Given $0<h<1$ and a known $0<G_0\le |D(B_0)|$, an initial
rank-one perturbation of norm at most $h/n^{10}$ makes
the new spectrum disjoint from the reference spectrum and retains at least
$1-1/(500n^3)$ of the discriminant modulus. For a subsequent pencil
$p-sq$, assume $p(\lambda_i)\ne0$ at every reference node.
If a parameter on the radius-$h$ circle attains at least half the
maximum discriminant modulus $M$, its radius can be increased by a
relative amount at most $2/(1000n^3)$ so that the discriminant modulus
decreases by at most $M/32$ and the characteristic polynomial remains
bounded away from zero at every reference eigenvalue.
\end{lemma}
\begin{proof}
Given $0<G_0\le|D(B_0)|$ and $0<h<1$, set
\[
 \alpha=\frac1{n(1+\|X_0\|_F^2)},\qquad
 \beta=\frac1{n(1+\|W_0\|_F^2)},\qquad
 u_0=\alpha X_0{\bf1},\qquad v_0^*=\beta{\bf1}^TW_0 .
\]
These vectors have norm at most one and can be computed in $O(n^2)$
operations. The residue at each reference eigenvalue is $\alpha\beta$. For
\[
 \sigma=\min\left\{h/n^{10},\,
           G_0/(2^{10}n^3\,8^{n(n-1)})\right\}
\]
the characteristic polynomial after adding $\sigma u_0v_0^*$ satisfies
$p_\sigma(\lambda_i)=-\sigma\alpha\beta p_0'(\lambda_i)\ne0$.
On the unit parameter circle the matrix norm is at most three, so its
discriminant modulus is at most $8^{n(n-1)}$. Cauchy's coefficient
bounds give a change of at most $2\sigma 8^{n(n-1)}$ at $\sigma$,
which proves the first assertion.

For a subsequent parameter $s_0$, consider $s_j=s_0(1+j/(1000n^4))$ for
$0\le j\le2n$. The discriminant polynomial has degree at most
$D=2n-2$. Its modulus on $|s|\le h(1+1/D)$ is at most $eM$ by the
exterior maximum principle. Discs of radius $h/(2D)$ about the
outward segment fit in this larger disc, giving a derivative bound
$6DM/h$ and a loss below $M/32$. For each reference eigenvalue, let
$a=p(\lambda_i)\ne0$, $b=s_0q(\lambda_i)$. We require
\[
 |a-b(1+j/(1000n^4))|\ge |a|/(10000n^4).
\]
If $|b|\le |a|/2$, every parameter satisfies this inequality.
Otherwise, consecutive values are separated by at least
$|a|/(2000n^4)$, so the inequality can fail for at most one index.
Binary search finds the allowed integers nearest to $1000n^4\Re(a/b-1)$.
At each such integer $j$, we test the inequality by comparing squared moduli.
Exclude the failing index for each reference eigenvalue and choose
any remaining index. After the multipoint evaluations, this
search costs $O(n\log n)$. An accepted step divides the previous
lower bound by at most $10000n^4$. When the accepted matrix also has norm at most two, the modulus of its monic characteristic
polynomial at $\lambda_i$ is a product of $n$ distances, each at most
four. Dividing its modulus bound by $4^{n-1}$ gives a lower bound
on every distance from a reference eigenvalue to the current spectrum.
\end{proof}

\subsection{Resolvent vectors and nearest gaps}
The reference resolvent gives approximations to the right and left
eigenvectors, together with the scalar product needed to estimate the
norm of a spectral projector. In the reference coordinates, the
perturbed matrix has the form
\[
 W_0BX_0=\Lambda+QFS^* .
\]
We choose $Q,S$ to have orthogonal columns of lengths in $[1/2,1]$.
Apply exact Gram--Schmidt to each list of perturbation factors in the
reference eigenbasis, omitting zero residuals and using the field
normalisation of Section~\ref{sec:foundations}. Coordinates in each
basis are inner products divided by squared column lengths. Absorbing
these coordinate maps into $F$ gives
$\|F\|\le4\|W_0(B-B_0)X_0\|$, even when the original columns are
linearly dependent. Let $r,s\le\min(k,n)$ denote the numbers of
columns of $Q,S$. Thus $Q\in\mathbb C^{n\times r}$,
$S\in\mathbb C^{n\times s}$ and $F\in\mathbb C^{r\times s}$. We assume below that the current spectrum is simple
and disjoint from the reference spectrum. In particular, both ranks
are positive.

For the accuracy estimates, fix a power of two $H$ bounding
$100,n^4$, the norms of $X_0,W_0,F$, and the reciprocals of the
current gaps and the distances to the reference spectrum.
Section~\ref{sec:global-assembly} gives such a bound, with $\log H$
polynomial in $n$ and $\log(n/\delta)$. For
$R(z)=(zI-\Lambda)^{-1}$, define
\[
 J(z)=S^*R(z)Q,\qquad L_R(z)=I_r-FJ(z),\qquad L_L(z)=I_s-J(z)F.
\]
Here $I_r$ and $I_s$ are the identity matrices of orders $r$ and $s$.
Away from the reference spectrum, both determinants equal
$\chi_B(z)/\chi_{B_0}(z)$. At an approximation $z_j$ matched to an eigenvalue of $B$,
choose a column $c$ of $\operatorname{adj} L_R(z_j)$ containing an
entry of largest modulus. Independently, choose a row $d$ of
$\operatorname{adj} L_L(z_j)$ by the same rule. In the original
coordinates, the corresponding column $x$, row $y$ and scalar pairing are
\begin{equation}
x=X_0R(z_j)Qc,\qquad y=dS^*R(z_j)W_0,\qquad
 yx=-dJ'(z_j)c .
\label{eq:resolvent-vectors}
\end{equation}
Here $y$ is a row, so $yx$ is ordinary row--column multiplication.

\begin{lemma}[Accuracy of the resolvent vectors]\label{lem:transfer-vectors}
Suppose $|z_j-\mu_j|\le H^{-60n}$ for a true eigenvalue $\mu_j$.
Evaluating the same column and row indices at $\mu_j$ gives nonzero
right and left eigenvectors. The relative errors in their norms and
in the pairing in \eqref{eq:resolvent-vectors} are each less than $1/(100n^2)$.
The Euclidean vector errors satisfy the same relative bound.
\end{lemma}
\begin{proof}
We first bound a cofactor from below at the true eigenvalue. On
$|z-\mu_j|\le1/(2H)$, the resolvent $R$, the two small matrices and
their derivatives are bounded by $2H,H^3,H^4$, respectively. At
$\mu_j$, the determinant derivative satisfies
\[
 |\chi'_B(\mu_j)/\chi_{B_0}(\mu_j)|\ge H^{-2n}
\]
using the lower bound on current gaps and the upper bound four on
distances to the reference spectrum. Jacobi's formula therefore gives
an adjugate entry of modulus at least $H^{-5n}$ in each small matrix.
Throughout the disc, the adjugates and their derivatives have norms at
most $H^{4n}$ and $H^{6n}$. To see the derivative bound, expand a
cofactor derivative as at most $n!n$ products, each containing at most
$n-2$ entries bounded by $H^3$ and one derivative bounded by $H^4$.
The scalar adjugates are constant. Consequently, the column and row
selected at $z_j$ each retain an entry of modulus at least
$H^{-5n}/2$ at $\mu_j$, including when the maximum is attained at
more than one entry.

At the true root, the adjugate identity gives $L_R c=0$ and
$dL_L=0$. These imply the eigenvector equations for the vectors in
\eqref{eq:resolvent-vectors}. The least
singular values of $Q,S$ on their column spaces are at least $1/2$,
and $\sigma_{\min}(R(\mu_j))\ge1/4$. Consequently each true vector
has norm at least $H^{-8n}$. On the disc, the vector functions with the selected indices have norms
at most $H^{6n}$ and derivatives at most $H^{8n}$.
For the true simple eigenvalue the spectral projector has norm at
most $(4H)^{n-1}$ by its polynomial product formula. The modulus of the pairing of the unit left and right vectors is its
reciprocal and is at least $H^{-2n}$, so the true unnormalised pairing
is at least $H^{-18n}$ in modulus.
The derivative of the product has magnitude at most $H^{17n}$.
Combining these estimates with the assumed root accuracy proves the
claims. If $c$ or $d$ is subsequently rescaled, apply the same
rescaling to the vector evaluated at the exact eigenvalue and to the
corresponding factor in the pairing.
\end{proof}

The values of $J,J'$ at all $n$ evaluation points require $O(k^2)$ rational
multipoint evaluations. For $J'$ we use the quotient rule, differentiating
the numerator and denominator. Direct cofactor determinants with exact pivoting cost
$O(nk^5)$, including when an evaluation point $z_j$ is an exact root. The orthogonal
factorisations and all operations on the small matrices therefore cost
$\widetilde O(nk^2)+O(nk^5)$. The scalar product in \eqref{eq:resolvent-vectors} is evaluated at
the same accuracy as the resolvent vectors.

\paragraph{Estimating nearest gaps.}
We estimate the Euclidean gaps by nearest-neighbour $\ell^1$
distances. For each of four reflected orientations, first consider
points in the northeast quadrant of an evaluation point. Sort by decreasing real
part, inserting each complete tie group before making the evaluations.
A segment tree over the sorted imaginary parts returns the two smallest
values of $u+v$ at distinct indices among points $(u,v)$ for which
$v$ is at least the evaluation point's imaginary part. Excluding the evaluation index
and subtracting the coordinate sum of the evaluation point gives the nearest distance in that
quadrant. Minimising over the four orientations gives the exact nearest
$\ell^1$ distance in $O(n\log n)$. This is between the Euclidean
gap and $\sqrt2$ times that gap. Moving every node by $H^{-60n}$
changes its $\ell^1$ gap by at most $2\sqrt2 H^{-60n}$.

\subsection{Changing the reference eigensystem}

If the approximate eigenbasis is invertible, its residual determines a
matrix correction that makes the diagonalisation exact. The following
proposition also constructs the inverse within the required cost.

\begin{proposition}[Explicit diagonalisation after a small correction]
\label{prop:batch-adoption}
Under the preceding assumptions on the disjoint reference and simple
current spectra, let $b\in\mathbb Z_{\ge0}$. One can compute a correction of
Frobenius norm at most $2^{-b}$ and an explicit diagonalisation whose
eigenvalues are field elements matched to the current eigenvalues.
It suffices to compute these approximations to accuracy
$2^{-b}H^{-Cn}$, for a sufficiently large absolute constant $C$.
After the scalar root approximation, the cost is
$\widetilde O(n^2(k+1)^{\omega-2}+nk^5)$.
\end{proposition}
\begin{proof}
Let $Z=\operatorname{diag}(z_j)$. We first form the approximate eigenbasis.
The matrix of unscaled coordinates $Y_{:j}=R(z_j)Qc_j$ has displacement
$\Lambda Y-YZ=-Q[c_j]$, so structured multiplication on the left
gives the columns $X_0Y$ in the original coordinates. For each column, component
rescaling and squared bisection give an upper norm estimate within a
factor $9/8$. Divide the column in the original coordinates, its column in the reference eigenbasis and its generator column by this estimate. Henceforth, $Y,X=X_0Y$ denote the rescaled matrices. The columns of $X$ have norms in $[8/9,1]$.
By Lemma~\ref{lem:transfer-vectors}, the true columns scaled by the
same divisors have norms in $[1/2,2]$. Before scaling, the residual
is exactly
\[
 (B-z_jI)x_j=-X_0Q\det L_R(z_j)e_{\ell(j)},
\]
where $\ell(j)$ is the selected column index. The residual is scaled
by the same divisor as the column. Hence
$\Delta=(XZ-BX)X^{-1}$ has a factorisation of small inner dimension. The Gram matrices of
these factors give its squared Frobenius norm.
Indeed, for $\Delta=UV^*$ this norm squared is
$\operatorname{tr}((U^*U)(V^*V))$. Forming $X_0Q$ and the other
products of dense matrices with these factors costs
$O(n^2(k+1)^{\omega-2})$ by blocking.

To bound the correction, consider the corresponding scaled eigenbasis
$X_*$. The $j$th dual row has norm at most $2\|P_j\|$, so
$\|X_*^{-1}\|\le2\sqrt n(4H)^{n-1}$. The derivative estimate and
the lower norm bound before scaling give
$\|X-X_*\|_F\le H^{O(n)}\max_j|z_j-\mu_j|$. Using
$BX_*=X_*\operatorname{diag}(\mu_j)$ gives the same estimate for
$\|XZ-BX\|_F$. The Neumann series for $(I+X_*^{-1}(X-X_*))^{-1}$ bounds $\|X^{-1}\|$ by
$2\|X_*^{-1}\|$ and proves the required correction bound. Distinctness of the eigenvalue approximations and the bound on the
correction norm can both be tested by field comparisons.
Once these tests hold, $B+\Delta=XZX^{-1}$ exactly.

It remains to construct the inverse within the stated cost. We first
introduce the real auxiliary nodes $a_i=8+i/n$, $b_i=12+i/n$.
Let $(H_1)_{ij}=1/(a_i-\lambda_j)$,
$(H_2)_{ij}=1/(z_i-b_j)$. The Cauchy determinant formula shows
that both matrices are invertible. For $Y'=H_1YH_2$, the identity
\[
 D_aY'-Y'D_b=(D_aH_1-H_1\Lambda)YH_2+
 H_1(\Lambda Y-YZ)H_2+H_1Y(ZH_2-H_2D_b)
\]
gives generators of width $k+2$, since the two outer displacements
have rank one. The factors in these generators are obtained by structured matrix
products.
Compute all entries of $Y'$ by a blocked generator product and
entrywise division by $a_i-b_j$. One structured multiplication with
disjoint node lists then gives $G=Y'^*Y'$, including its diagonal.
Since $a,b$ are real, conjugate transposition and the product rule give
the displacement $D_bG-GD_b$, with the same diagonal matrix on both
sides, of width
$w=2(k+2)$.

Suppose that a matrix of order $q$ has a known displacement with the
same diagonal matrix on both sides, and that its diagonal is also
known. To multiply this matrix, divide its nodes into two halves. The off-diagonal blocks
have disjoint node lists, while the diagonal blocks inherit their
diagonals. Repeat this subdivision. Multiplication with $q$ dense
right-hand sides costs $\widetilde O(q^2w^{\omega-2})$: the block
row counts at each level sum to $q$, and subdivision stops when a block first has size at most $w$.
At that level square blocking forms the off-diagonal entries from the
generators, using the known diagonal on blocks with identical node lists, and multiplies
the resulting dense blocks by the right-hand sides.

To invert $G$, partition a current block as
$\left(\begin{smallmatrix}M&B'\\ C'&D'\end{smallmatrix}\right)$
and partition its two displacement generators as $U=(U_1;U_2)$ and
$V^*=(V_1^*,V_2^*)$. We must retain the diagonal of the Schur
complement as well as its generators, since this displacement does not
determine diagonal entries. Compute it from
$D'-C'M^{-1}B'$. Once $C'M^{-1}$ is available, the diagonal
correction requires row--column inner products and costs $O(q^2)$
operations. The Schur complement has generator factors
\[
 U_2-C'M^{-1}U_1,\qquad V_2^*-V_1^*M^{-1}B'.
\]
To verify this formula, multiply the full displacement on the left by
$[-C'M^{-1},I]$ and on the right by $\binom{-M^{-1}B'}{I}$.
The generator width passed to each recursive inverse is therefore
still $w$. Once a block has been explicitly inverted, multiplying its
old generators by the inverse on both sides and negating one factor
gives generators of the inverse. By the product rule, the
off-diagonal products $C'M^{-1}$ and $M^{-1}B'$ have disjoint
endpoint node lists and width at most $2w$. Compute their entries by
multiplying the structured off-diagonal factors by the dense inverse.
These representations give the Schur complement and the four standard
inverse blocks, including
$M^{-1}+(M^{-1}B')\,(D'-C'M^{-1}B')^{-1}\,(C'M^{-1})$.
The intermediate generator widths are at most $3w$ by the product
rule. Retain the diagonals of the explicitly formed inverse blocks
for subsequent products involving identical node lists. Stop the inverse recursion when the block order is at most $w$, and use
dense inversion there. At the last split the two child orders are
comparable to $w$, so their entries can be recovered from the generators
and diagonal by square blocking. At larger orders, each multiplication of a dense
matrix by a generator of width $w$ costs $O(q^2w^{\omega-2})$ for $q\ge w$
by blocking. The remaining products at a recursion node cost
$\widetilde O(q^2w^{\omega-2})$ by the preceding subdivision.
Summing the squared block sizes over the recursion gives a geometric
series. Every required principal pivot exists because the Gram matrix
is positive definite. If the initial order is at most $w$, the Gram matrix is already explicit
and a dense inverse costs $O(n^\omega)$. Dense inversion of positive
definite blocks has this cost by recursive Schur complementation with
square multiplication. A failed pivot in the Gram
recursion detects a singular approximate eigenbasis.

Finally $Y^{-1}=H_2G^{-1}Y'^*H_1$ is formed with the structured
multipliers. Multiplying $\Lambda Y-YZ$ on both sides by $Y^{-1}$
and negating gives generators for $ZY^{-1}-Y^{-1}\Lambda$.
A final structured multiplication with disjoint node lists yields $X^{-1}=Y^{-1}W_0$.
\end{proof}

\section{Norms in the original coordinates}
\label{sec:mixed-range-bridge}

Since the reference eigenvectors need not be orthogonal, Euclidean norms
in the reference coordinates do not determine the norms of the
corresponding eigenvectors. Multiplication by the reference basis for
every eigenvector after each update would be too expensive. Instead,
we estimate all the norms with an accuracy adapted to each reference
column, and compute the vectors in the original coordinates only for
selected eigenvalues.

\subsection{Reference precisions and integer packing}

Suppose the reference matrix $B_0=X\Lambda X^{-1}$ has $|\lambda_i|<2$,
$1/2\le\|x_i\|\le2$, and let $w_i^*$ denote the dual rows of $X$.
Here $X$ is the reference basis denoted by $X_0$ in
Section~\ref{sec:transfer-rebase}, and $W_0=X^{-1}$.
Let $s_i^{(0)}=\min_{j\ne i}|\lambda_i-\lambda_j|$. We use positive
integer precision exponents $b_i$, each a power of two, so that the
corresponding error scales are $2^{-b_i}$. Assume that
\[
 b_i\ge20\log_2(2n)+\log_2^+\|w_i\|+\log_2^+(1/s_i^{(0)}).
\]
Successive squaring and bracketing by squared comparisons give these
exponents within a constant factor of the right-hand side. Throughout
this section, assume $\sum_i b_i\le CnL$, where
$L\ge\log_2(2n/\delta)$ and $C$ is absolute. Increasing the $b_i$
by an absolute factor gives the same hypothesis for normalised left
eigenvectors, since the dual norm on either side is comparable with the
same projector norm. We require this sum bound only for the reference
eigensystem.
Lemma~\ref{lem:total-precision} verifies this hypothesis whenever the
local algorithm starts a batch at a reference basis.

At a fixed precision we can scale the rounded entries of a matrix and a
vector to bounded integers. Under the unit-cost convention we can encode each matrix column as a
single integer. The digits of a linear combination of these integers
then give the entries of the product, up to known additive corrections.

\begin{lemma}[Integer packing]\label{lem:packed-matvec}
Let $H\ge1$ be an integer, and let $M$ be a fixed $a$-by-$b$
integer matrix with $a,b\le n$ and entries in $[-H,H]$. After
preprocessing its columns, a product with an integer vector
$v\in[-H,H]^b$ can be computed in $O((a+b)(1+\log D))$ field
operations and comparisons, for any power of two $D>4nH^2$.
\end{lemma}
\begin{proof}
Number the rows from zero, and retain
$C_j=\sum_i(M_{ij}+H)D^i$ together with the row sums. The digit
in position $i$ in the base-$D$ expansion of
$\sum_j C_j(v_j+H)$ is $\sum_j(M_{ij}+H)(v_j+H)$, since this sum
lies in $[0,D)$ and there is no carry. Starting at the highest position,
find each digit by binary search and subtract its multiple of $D^i$.
The row sums and $\sum_jv_j$ remove the biases. Preprocessing takes
$O(a)$ operations to tabulate powers of $D$ and $O(ab)$ to insert
already quantised entries.

For real fixed-point inputs with $O(p)$ magnitude and precision bits,
scale by powers of two and quantise each new scalar using $O(p)$
squared or ordinary comparisons. If $p\ge\log_2(2n)$, the cost of
each product after preprocessing is $O((a+b)p)$. Four real products
suffice for complex pairs.
\end{proof}

\subsection{Sign directions and Frobenius inner products}

\begin{lemma}[Sign directions in the reference basis]\label{lem:batch-signs}
Under the hypotheses on $X$ and $b_i$ above, set
$\epsilon=2^{-40}n^{-6}(1+L)^{-2}$. After $\widetilde O(n^2 L)$
preprocessing, the following construction costs
$\widetilde O(nm+nL)$ operations. Given approximations to $m$ pairs of unit left and right eigenvectors
of the current matrix, with the accuracy permitted by
Lemma~\ref{lem:signs},
one obtains sign vectors $\bar u,\bar v$ with entries $\pm1/(2s)$
as in that lemma and field vectors $\mathbf a,\mathbf b$ in reference coordinates.
The directions $u=X\mathbf a$ and $v^*=\mathbf b^*X^{-1}$ satisfy
\[
 \|u-\bar u\|,\ \|v-\bar v\|\le\epsilon,\qquad \|u\|,\ \|v\|\le1 .
\]
Here $P_i$ is the spectral projector of the current matrix for a
selected eigenvalue. For at least $m/25$ of these eigenvalues, $|v^*P_i u|\ge\|P_i\|/(2^{13}n)$.
\end{lemma}
\begin{proof}
Let $W=X^{-1}$, and denote its rows by $w_i^*$. Choose a dyadic
$e\in[\epsilon/(32n^2),\epsilon/(16n^2)]$. Quantise the real
and imaginary parts in row $i$ of $W$ on the mesh $e$, and in
column $i$ of $X$ on the mesh $e\,2^{-b_i}$. Denote the resulting
arrays by $\widehat W,\widehat X$. Since $\|w_i\|\le2^{b_i}$ and
$\|x_i\|\le2$, the integer entries after scaling in group
$b_i=b$ require $O(b+\log(1/e))$ bits. Entrywise binary search
and column packing therefore cost
$O(n\sum_i(b_i+\log(1/e)))=\widetilde O(n^2L)$.

Apply Lemma~\ref{lem:signs} to the approximate eigenvectors. Take
$\mathbf a=\widehat W\bar u$ and $\mathbf b^*=\bar v^*\widehat X$. To obtain either
product, group the output entries by $b_i$ and perform the sign
product before division by $2s$. In each real integer array group
use the packed columns $C_j$ of Lemma~\ref{lem:packed-matvec}.
With unscaled input signs $\sigma_j$, form
$\sum_j C_j(1+\sigma_j)$ in $O(n)$ operations. The digits determine
the signed row products after subtraction of the known row sums and
biases. A radix exceeding $4n(H+1)^2$ suffices for entries bounded by
$H$. Recovering one digit costs $O(b+\log(1/e))$ comparisons.
There are $O(\log(2+nL))$ groups, since the $b_i$ are powers of two.
The total cost of the products after choosing the signs is $\widetilde O(nL)$.

In position $i$ the errors in $\mathbf a$ and in the row $\mathbf b^*$ are at most
$2ne$ and $2ne\,2^{-b_i}$, respectively.
Summing the first errors against the columns of $X$ and the second
against the rows of $W$ gives the claimed Euclidean errors. Since the original sign vectors have norm at most $1/2$, the
constructed directions have norm at most one. At the indices where the
sign construction succeeds, their overlaps with the unit eigenvectors have
modulus at least $1/(16s)-\epsilon\ge1/(32s)$. Their product times
$\|P_i\|$ is $|v^*P_i u|$, and $s<2\sqrt n$. The two reference
factors $\mathbf a,\mathbf b^*$ define the rank-one update exactly, so its pencil
can be formed without multiplying either factor by a dense basis.
\end{proof}

\begin{lemma}[Frobenius inner products in a batch]\label{lem:proxy-query}
Let $F_0$ be an explicit matrix with $\|F_0\|_F\le\delta$.
During a batch let
$F=F_0+\sum_{j=1}^l t_j\bar u_j\bar v_j^*$, $l\le k\le n$,
with sign directions as in Lemma~\ref{lem:batch-signs}. For new
sign directions $\bar u,\bar v$ and a real $c$ satisfying
$|c|\delta\le 2^{10}n^2(1+L)$, the number
$c\,\operatorname{tr}(F^*\bar u\bar v^*)$ can be obtained to
absolute error $1/100$ in $\widetilde O(nL+nk)$ operations after
$\widetilde O(n^2 L)$ preprocessing. At the end of the batch, forming the matrix $F$ explicitly costs $O(n^2 k^{\omega-2})$ operations.
\end{lemma}
\begin{proof}
Quantise the real and imaginary parts of $F_0/\delta$ to within
$2^{-p}$, where $p=40+10\lceil\log_2(2n(1+L))\rceil$. Pack its
columns as in Lemma~\ref{lem:packed-matvec}. Apply this matrix to
the unscaled signs for $\bar v$, divide by $2s$ and contract with
$\bar u$. After conjugation and multiplication by $c\delta$ this
has error at most $2|c|\delta n2^{-p}<1/100$. The contribution of
each stored term is exactly
$c\overline{t_j}(\bar u_j^*\bar u)(\bar v^*\bar v_j)$.
Finally multiply the two $n$-by-$l$ arrays of directions with the
scalars $t_j$ absorbed into the first array. Blocking the two outer
dimensions in intervals of length $l$ gives the asserted cost.
When $l=0$ the initial matrix already gives the result.
\end{proof}

\subsection{Triangular factorisation with individual precisions}

For the triangular factorisation we take the columns in decreasing order
of the precision they require. Two passes of orthogonalisation control
both the loss of orthogonality and the error in reconstructing each column.

\begin{lemma}[A triangular factorisation with individual accuracies]
\label{lem:mixed-qr}
Permute the reference columns and dual rows so that
$b_1\ge\cdots\ge b_n$. Apply the same permutation to the reference
eigenvalues and to the rows of every numerator factor in these
coordinates. The eigenvalue approximations retain their original ordering. Apply
the construction separately to the right eigenvectors and the conjugate
transposes of the left eigenvectors, using a consistent permutation of
the reference data on each side.
In $O(n\sum_i b_i)$ operations up to fixed logarithmic factors one
can compute $Q,R\in\mathbb C^{n\times n}$, with $R$ upper triangular, such that
\begin{align}
 |(Q^*Q-I)_{ij}|&\le2^{-2b_j} &&(i\le j),\label{eq:mx-Q}\\
 \|x_j-QR_{:j}\|&\le2^{-8b_j}.\label{eq:mx-X}
\end{align}
In particular, for every $y$,
\begin{equation}
\|QRy-Xy\|\le
 \Big(\sum_i 2^{-8b_i}\|w_i\|\Big)\|Xy\|,\qquad
 \|Q^*Q-I\|\le n2^{-2\min_i b_i}.
\label{eq:mixed-reconstruction}
\end{equation}
\end{lemma}
\begin{proof}
Let $q_j$ be column $j$ of $Q$ and
$Q_{<j}=[q_1,\ldots,q_{j-1}]$. For column $j$, apply two
orthogonalisation passes to $x_j$ against $Q_{<j}$. On the current residual $u$, each pass computes
$a\simeq Q_{<j}^*u$ and subtracts a packed approximation to
$Q_{<j}a$. Absolute component precision $2^{-30b_j}$ on bounded
inputs gives Euclidean error below $2^{-20b_j}$ in each product,
including the dimension factors. The sum of the two computed vectors
$a$ gives the strictly upper part of the current column of $R$.

We first show that the residual remains bounded away from zero.
Inductively, let $T$ map $x_i$ to $QR_{:i}$ for $i<j$ and fix the
other basis columns. The reconstruction bound and duality give
$\|T-I\|\le\sum_{i<j}2^{-8b_i}\|w_i\|\le n2^{-7b_j}$.
Since $w_j^*$ annihilates the preceding columns and $w_j^*x_j=1$,
the distance of $x_j$ to their span is at least $1/\|w_j\|$.
Their transformed span equals the span of $Q_{<j}$ by the nonzero
triangular diagonals. As $T$ fixes $x_j$, its distance to this
span is at least $(1-n2^{-7b_j})/\|w_j\|$. The errors in the two vector subtractions together have norm below $2^{-18b_j}$, so the final residual
has norm at least $0.8\,2^{-b_j}$.

Let $e=\|Q_{<j}^*Q_{<j}-I\|\le n2^{-2b_j}$. Each exact pass
multiplies the vector of inner products by $I-Q_{<j}^*Q_{<j}$, and
$\|I-Q_{<j}Q_{<j}^*\|\le1$. Since $\|x_j\|\le2$, these bounds
keep the residual and coefficient norms below four, as required by the
packing estimates. After two passes,
$\|Q_{<j}^*u\|$ is at most $4e^2+O(2^{-20b_j})$. Squared
comparisons compute $d_j>0$ approximating $\|u\|$ to relative error
$2^{-10b_j}$ in $O(b_j)$ operations. Set $q_j=u/d_j$,
$R_{jj}=d_j$. Division by $d_j$ makes the inner-product error less
than $2^{-2b_j}$, and the squared column norm has the same error
bound. Since $R$ uses the computed coefficients $a$, only the errors in the two vector subtractions enter the reconstruction error. This proves
\eqref{eq:mx-Q} and \eqref{eq:mx-X}. The first bound in \eqref{eq:mixed-reconstruction} follows
from $y_i=w_i^*Xy$ by summing the column errors.

For the cost bound, retain packed columns for $Q_{<j}$ and
$Q_{<j}^*$ at every required dyadic precision up to $b_j$, using a
fixed radix for each precision. Appending $q_j$ adds one column to
the first table and one digit to each of $n$ columns in the second.
The active dyadic precisions sum to $O(b_j)$. Thus each step,
including quantisation, costs $O(nb_j)$ up to logarithms.
\end{proof}

\subsection{Euclidean norms and selected directions}

We apply this factorisation at $n$ distinct evaluation points $|z|\le2$
that do not belong to $\{\lambda_i\}$. At each point, the coordinate
vector and its norm in the original coordinates have the form
\begin{equation}
y_i(z)=\frac{(A c)_i}{z-\lambda_i},\qquad \mathcal N=\|Xy(z)\|.
\label{eq:coordinate-vector}
\end{equation}
Here $A$ has $1\le r\le n$ columns, and $c$ is a field vector
depending on $z$. For the right vector in \eqref{eq:resolvent-vectors}, $A=Q$ is the
resolvent factor. In the rest of this section, however, $Q,R$ denote
the triangular factors of Lemma~\ref{lem:mixed-qr}.

For the left row in \eqref{eq:resolvent-vectors}, apply the construction
to its conjugate transpose. The evaluation points as well as the reference nodes
are then conjugated, and the selected cofactor row is replaced by its
conjugate transpose. Let $D$ be the positive diagonal matrix that
normalises the columns of $W_0^*$ to within the prescribed constant
factor. Its basis in the original coordinates is $W_0^*D$, and the
numerator factor in that basis is $D^{-1}S$. In both cases, all data
in \eqref{eq:coordinate-vector}, including the computed evaluation point $z$, are exact field values.

We compute the triangular coordinates for the rows requiring higher precision
exactly. For the remaining rows we estimate their contribution to the
squared norm from moments of the Gram matrix.

\begin{proposition}[Estimating all Euclidean norms]\label{prop:mixed-norms}
Let $\tau\ge1$ be a power of two and let $d=\#\{i:b_i>\tau\}$. After
$\widetilde O(n\sum_i b_i+n^2(\tau+\log n))$ preprocessing for the reference eigensystem,
the $n$ squared norms in \eqref{eq:coordinate-vector} can be estimated to relative error
$1/20$ in
\[
 \widetilde O\big(n\max(d,r)r^{\omega-2}
            +n(r+1)^2(\tau+\log n)^5\big)
\]
field operations. Here and below logarithmic factors include powers
of $1+\log(n/\delta)$ wherever $L$ occurs in the application.
\end{proposition}
\begin{proof}
The sum bound gives $d=O(nL/\tau)$. We first compute the initial
$d$ coordinates of $Ry(z)$ exactly at every evaluation point. The matrix
whose columns are the $y(z)$ has Sylvester displacement $-A[c]$
with respect to the reference nodes and evaluation points. Transposed structured
multiplication from Section~\ref{sec:transfer-rebase} gives the first
term in the cost. By upper triangularity, reference columns with
$i\le d$ contribute nothing below the first $d$ output coordinates.

For the remaining coordinates, let $R_t$ be the submatrix with row
and column indices greater than $d$. Lemma~\ref{lem:mixed-qr} gives
$\|R_{:i}\|\le(\|x_i\|+2^{-8b_i})/\sigma_{\min}(Q)<3$.
Thus the entries of $R_t^*R_t$ have modulus at most nine and
$\|R_t\|\le 3\sqrt n$. Quantise $R_t$ using
$p+20\lceil\log_2(2n)\rceil$ fractional bits and apply its packed
adjoint to each column. This approximates $G=R_t^*R_t$ to absolute
entry error $2^{-p}$ in $O(n^2p)$ operations, where
$p=C_1(\tau+\lceil\log_2(2n)\rceil+64)$ and the absolute constant
will be chosen below. This Gram matrix depends only on the reference
basis.

We also need a bounded representation of the numerator on the remaining
rows $i>d$. Let $r_0=\operatorname{rank}(A_{\rm low})\le r$.
Exact Gram--Schmidt, omitting dependent columns, gives
$A_{\rm low}=UH_A$, with $U\in\mathbb C^{(n-d)\times r_0}$
and $H_A\in\mathbb C^{r_0\times r}$. Retain the original reference
indices $i>d$ for the rows of $U$. Before bracketing
the norm of a nonzero residual, scale it by its largest real or
imaginary component in modulus. Thus the columns of $U$ can be
scaled to have lengths in $[1/2,2]$ in $O(nr^2+r\log n)$
operations, independently of how small the singular values of
$A_{\rm low}$ are. For $c'=H_Ac$, we obtain
\begin{equation}
|y_i|\le\|w_i\|\mathcal N\le2^\tau \mathcal N,\qquad
 \|c'\|\le2\|A_{\rm low}c\|\le8\sqrt n\,2^\tau \mathcal N
 \quad (i>d),
\label{eq:coordinate-bounds}
\end{equation}
using $|z-\lambda_i|\le4$. Empty coordinate sets and rank-zero
numerators contribute zero.

We now group reference nodes so that their reciprocal denominators admit
short expansions. Cover the node disc by a fixed square. Subdivide
each square containing remaining reference nodes into four congruent squares until the
side length $w$ lies between $2^{-\tau}/4000$ and
$2^{-\tau}/1000$. The resulting tree has $T=O(\tau+\log n)$
levels. Assign boundary points to a unique child by a fixed half-open
convention (including the exterior edges of the initial square).
Midpoint comparisons assign each reference node to its path in $O(T)$ operations.
Store the occupied cells at each level by their integer row and column
indices, obtained by doubling the parent indices and adding the child
offsets. Sorting these lists permits neighbour searches by binary search
in $O(\log(2n))$ comparisons. For an evaluation point, accept a cell $I$ with side $w_I$ and centre
$a_I$ if $|z-a_I|>8w_I$, and subdivide it otherwise. Only a bounded
number of cells are near the evaluation point at each level, so only a bounded
number of their children can be accepted at the next level. At the
finest level, all remaining near points are within $9w$ of the evaluation point.
Since they are separated by at least $2^{-\tau}$, there is at most
one. Use its coordinate unit vector and compute its scalar $y_i$
exactly.

On an accepted cell, approximate the reciprocal denominator by
\begin{equation}
\sum_{j=0}^{s-1}
 \left(\frac{\lambda_i-a_I}{w_I}\right)^j
 \frac{(w_I/(z-a_I))^j}{z-a_I},\qquad i\in I .
\label{eq:reciprocal-expansion}
\end{equation}
The relative error is $|(\lambda_i-a_I)/(z-a_I)|^s
\le8^{-s}$. Hence the error in the remaining coordinate vector is at most
$8^{-s}\sqrt n\,2^\tau \mathcal N$. Taking
$s=128(\tau+\lceil\log_2(2n)\rceil+20)$ makes the error after multiplication by $R_t$ less than $\mathcal N/1000$. The estimate
also holds when cancellation makes the exact remaining component zero.

The expansions give a small set of vectors whose Gram moments determine
the norm of the remaining component. For a cell $I$, use auxiliary vectors with entries
$U_{ih}((\lambda_i-a_I)/w_I)^j$ on $I$, for $1\le h\le r_0$ and
$0\le j<s$, and zero elsewhere. Include the finest-level singleton
coordinate unit vectors. Each evaluation point uses $O((r+1)(s+1)(T+1))$ auxiliary vectors.
We need the moments $f_J^*Gf_I$ only for pairs that can occur at a
common evaluation point.

Suppose $I$ is the coarser cell. The distance inequalities at the parent cells
place both centres within $O(w_I)$ of the evaluation point. Thus $J$ lies in
the union $N(I)$ of a fixed number of neighbouring grid cells at the
level of $I$. These neighbourhoods have bounded overlap within each
level. Retain the rounded submatrices $G_{N(I),I}$ in packed form.
Since their entries are already quantised, insertion costs $O(n^2 T)$
for the entire tree, once per reference basis.

For each auxiliary vector on $I$, a packed product gives $G_{N(I),I}f_I$.
At each finer or equal level, group the reference nodes in $N(I)$ by their
cell $J$ and compute the inner products with each type of auxiliary vector.
At each coarser level, the total number of reference nodes examined over all $I$ is $O(n)$.
Keeping the indices of the occupied cells sorted at each level, we obtain the cost
\[
 \widetilde O\big(n(r+1)^2(s+1)^2(T+1)^2(p+\log n)\big).
\]
The tree and the packed Gram submatrices are reused for every update
at the reference basis. Only the auxiliary vectors involving $U$ must be
recomputed after an update.

It remains to control rounding in these moments. Estimate the squared norm of the remaining component by the real part
of the quadratic form in the coefficients of the auxiliary vectors. Taking the real part cannot increase the error. By \eqref{eq:coordinate-bounds},
the minimum leaf size and $|w_I/(z-a_I)|<1/8$, every coefficient associated with an accepted cell
has modulus at most $2^{2\tau+16}n\mathcal N$. Each coefficient associated with a nearby node is exactly $y_i$
and is bounded by \eqref{eq:coordinate-bounds}. The entries of all auxiliary vectors have modulus at most two.

Quantising the bounded auxiliary vectors, Gram entries and packed-product
inputs with $p+20\lceil\log_2(2n)\rceil$ fractional bits bounds each
moment error by $128n^2 2^{-p}$. Each inner product defining a moment has at most $n^2$
terms, and the additional $20\lceil\log_2(2n)\rceil$ bits absorb
the errors in both auxiliary vectors. There are at most
$2^{14}n(s+1)(T+1)$ coefficients. Hence the quadratic-form error is
at most
\[
 2^{4\tau+67}n^6(s+1)^2(T+1)^2 2^{-p}\mathcal N^2.
\]
Taking $C_1=4096$ makes this less than $\mathcal N^2/1000$, even in the
presence of cancellation. Combine this with the truncation error,
the exact first $d$ coordinates and the near-isometry and
reconstruction bounds in \eqref{eq:mixed-reconstruction}. The resulting absolute error in the
squared Euclidean norm is at most $\mathcal N^2/20$. Evaluating the $n$
quadratic forms is covered by the moment construction bound. If there
are no remaining coordinates, their contribution is zero. If $d=0$, there
are no initial coordinates to compute exactly.
\end{proof}

Once the evaluation points have been selected, the sign construction also needs
their left and right directions in the original coordinates.

\begin{proposition}[Recovering selected directions]\label{prop:mixed-directions}
For the unscaled vector approximations in \eqref{eq:resolvent-vectors} represented by
\eqref{eq:coordinate-vector}, under the accuracy hypothesis of
Lemma~\ref{lem:transfer-vectors}, choose any
$1\le m\le n$ evaluation points with $\mathcal N>0$.
After preprocessing the reference eigensystem, one can compute
normalising factors $\mathcal N'\in[\mathcal N/2,2\mathcal N]$ and
approximations to $Xy/\mathcal N'$ with $\ell^2$ errors less than $n^{-3}$. In addition to the
norm estimates and $\widetilde O(nr^2+r\sum_i b_i)$ operations
for the numerator matrix $A$, the cost is
\[
 \widetilde O\!\left(m\Big(\sum_i b_i+
                 r^2\log(2+\max b_i)+n\log n+\log(2+n\log H)\Big)\right).
\]
\end{proposition}
\begin{proof}
Proposition~\ref{prop:mixed-norms} gives $\mathcal N'$ by squared comparisons
with dyadic scales. For these unscaled vector approximations,
Lemma~\ref{lem:transfer-vectors} bounds both $\mathcal N$ and $1/\mathcal N$ by
$H^{O(n)}$. Compare the norm estimate to squared dyadic scales, inverting the estimate
first if it is below one. Successive squaring brackets the binary exponent.
Using the stored powers in descending order determines that exponent in
$O(\log(2+n\log H))$ operations. A dyadic scale nearest to the estimated
length on a logarithmic scale is within a factor $\sqrt{2}$ of that
length, and hence within a factor two of $\mathcal N$. The comparison at the
geometric midpoint uses the squared norm estimate directly.

Group the indices $J$ with common $b_i=b$. As for the rows requiring less precision,
form $A_J=U_J H_J$ with rank $r_J$ at most $\min(r,|J|)$.
The factorisations cost $O(nr^2+r\log n\,\#\{J\})$, and quantising
and packing the matrices $U_J$ costs $O(r\sum_i b_i)$ operations.

At each selected evaluation point, compute $c_J=H_J c$ in $O(rr_J)$
operations. The argument for \eqref{eq:coordinate-bounds} gives
$\|c_J\|/\mathcal N'\le16\sqrt n\,2^b$. Scanning the group nodes finds any
pole within $2^{-b}/4$ of $z$. There is at most one, since two
group nodes are separated by at least $2^{-b}$. Compute $y_i/\mathcal N'$
at this pole exactly from \eqref{eq:coordinate-vector}. At the other nodes, the reciprocal
denominator is at most $4\cdot 2^b$.

A product $U_J(c_J/\mathcal N')$ at $O(b)$ fixed-point precision, followed
by exact divisions, therefore gives absolute error at most $n^{-12}$
in every normalised coordinate. The intermediate magnitudes and the
reciprocal denominators each require $O(b+\log n)$ bits, while
$b\ge20\log_2(2n)$. Also, the true normalised coordinates have modulus at
most $2^{b+1}$, since $|y_i|\le\|w_i\|\mathcal N$. The packed product costs
$O((|J|+r_J)b)$ per group, and the exact row at a near pole costs
$O(r)$. There are at most $O(\log(2+\max_i b_i))$ groups.

Let $d_b=\#\{i:b_i\ge b\}$. The group columns of $R$ vanish below
row $d_b$, so form their contribution to $Ry/\mathcal N'$ by a packed
matrix--vector product at $O(b)$ bits. The coefficients are bounded
and the input has range $2^{b+1}$. With sufficiently large absolute
precision constants, the total coordinate error in each block,
including the input errors, is at most $n^{-9}$. The work is
$O((|J|+d_b)b)$. Pack these blocks once per reference basis at a
cost of $O(n\sum_i b_i)$ operations. Since the precisions are dyadic,
$\sum_b d_b b\le2\sum_i b_i$.

Sum the block contributions. By \eqref{eq:mixed-reconstruction}, $\|Ry/\mathcal N'\|$ is bounded by an
absolute constant. One final packed product by $Q$ at $O(\log n)$
bits therefore suffices. Increasing the absolute precision constants
makes the accumulated Euclidean error, including the error in \eqref{eq:mixed-reconstruction},
less than $n^{-3}$. To obtain any prescribed fixed constant in the error bound, increase also
the lower bound on the $b_i$ by an absolute additive constant before
dyadic rounding. This controls the reconstruction error in
\eqref{eq:mixed-reconstruction} with that constant and preserves
$\sum_i b_i=O(nL)$.
\end{proof}

\subsection{Projector estimates and cost per update}

Let $P_i$ and $s_i$ denote the spectral projector and nearest gap at
the $i$th matched eigenvalue of the current matrix. To estimate
$\|P_i\|/s_i$, take the computed column $x_i$ and row $y_i$ in
\eqref{eq:resolvent-vectors} at this evaluation point. Let $N_{R,i},N_{L,i}$
be their squared norm estimates and $d_i>0$ the distance from the $i$th
computed node to its nearest neighbour in $\ell^1$, obtained as in
Section~\ref{sec:transfer-rebase}. With the same column and row scales
in the norm estimates and the scalar product computed without rounding, define
\[
 \widehat K_i^2=
 \frac{N_{R,i}N_{L,i}}{|y_i x_i|^2d_i^2}.
\]
By Lemma~\ref{lem:transfer-vectors} and Proposition~\ref{prop:mixed-norms},
$\widehat K_i$ is within a factor four of $\|P_i\|/s_i$. Real-field
operations give the squared modulus of the pairing, and all comparisons
can be made with $\widehat K_i^2$.

For selected evaluation points, Proposition~\ref{prop:mixed-directions} gives
vectors whose norms, before approximation error, lie between fixed
positive constants. Use the permitted accuracy $10^{-6}n^{-3}$ and
normalise by bisection on a bounded interval. Lemma~\ref{lem:transfer-vectors},
with the same accuracy for the eigenvalue approximations or a larger
absolute exponent, together with the individual approximation bounds,
gives unit directions within
$1/(100n^2)$ of the true directions. Lemma~\ref{lem:batch-signs}
then gives two sign vectors and exact reference coordinates for nearby
perturbation directions, at the stated cost.

For the iteration, suppose $r=O(n^b)$ and choose $\tau$ within a fixed
factor of $n^aL$, where $0<a<1$ and $0<b\le1/2$ are fixed. The norm
estimates cost
\[
 \widetilde O\bigl(n^{2-a+(\omega-2)b}
                 +n^{1+(\omega-1)b}+n^{1+2b+5a}\bigr)
\]
per update, after $\widetilde O(n^{2+a})$ preparation for the reference
basis. In the application $\log H=(256nL)^8\log 2$, so the additional scale
search costs $O(\log(nL))$ per evaluation point. Since $\sum_i b_i=O(nL)$ and
$r^2\le O(n)$, the cost of recovering
any $m$ selected directions is $\widetilde O(nm)$, in addition to
$\widetilde O(nr^2+rnL)$ for preprocessing the numerator matrix. No relation between
$m$ and $r$ is required. We use
\[
 a=\frac{23}{192},\qquad b=\frac{25}{144},\qquad
 \rho=\frac{545}{576}.
\]
The matrix multiplication bound $\omega<2.38$ from
\cite{alman-multiplication} gives
\[
 2-a+\frac{19}{50}b=1+2b+5a=2+a-b=1+\rho.
\]
Thus the norm estimates cost $\widetilde O(n^{1+\rho})$ per update,
and preparation divided over $k=\lceil n^b\rceil$ updates has the same
bound.

\section{Iteration and proof of the regularisation theorem}
\label{sec:global-assembly}

At each stage we compare each spectral projector norm with the
distance from its eigenvalue to the rest of the spectrum. When
sufficiently few ratios exceed the local bound, we retain the invariant
subspace spanned by the other right eigenvectors. Only the trailing
diagonal block then requires a recursive call.

Throughout this section, $n$ and $\delta$ denote the dimension and
tolerance of the current recursive call. We reserve $N$ and $\delta_*$
for the original input when deriving the final bounds.
In each call the input is a contraction of order $n$, with
perturbation tolerance $0<\delta<1/2$. The top-level call has depth zero,
and each recursive call increases the depth by one. A call in dimension
zero returns the empty matrix and basis. A call in dimension one returns
the given scalar and the unit column, with its inverse.
For $n\ge2$, a call is accelerated if its depth is zero or one
and $n\ge2^{20}$. At these two depths, respectively, set
\[
  \alpha=\begin{cases}947/1000,&\text{at depth zero},\\39/50,&\text{at depth one},\end{cases}
 \qquad m=\lfloor n^\alpha\rfloor .
\]
All other calls are terminal, with $m=1$. Define
\begin{equation}
 \begin{gathered}
 L=\lceil\log_2(64n/\delta)\rceil,\qquad d=n(n-1),\qquad
 T=2^{16}n^2L,\\
 \nu=2^{24}n^2L/\delta^2,\qquad
 \eta=\delta/(2^{30}n^2L),\\
 G_0=(\delta/(8n))^d/16,\qquad H=2^{(256nL)^8}.
 \end{gathered}\label{eq:parameters}
\end{equation}
Choose the largest dyadic radius $r=2^j$, $j\in\mathbb Z$, with
\begin{equation}\label{eq:local-radius}
 r^2\le 2^{-40}\frac{\delta^2m}{n^2L},
\end{equation}
and set
\begin{equation}\label{eq:local-thresholds}
 M=\frac{2^{20}n}{r},\qquad Q=32M,\qquad
 Q_{\rm term}=\frac{2^{12}}r,\qquad
 k=\lceil n^{25/144}\rceil .
\end{equation}
Thus $0<r<1$, and
\begin{equation}\label{eq:local-threshold-bounds}
 Q\le 2^{70}\frac{n^2\sqrt L}{\delta\sqrt m}
 \quad\hbox{in an accelerated call},\qquad
 Q_{\rm term}\le 2^{70}\frac{n\sqrt L}{\delta}
 \quad\hbox{in a terminal call}.
\end{equation}
The algorithm chooses $r$ by squared comparisons. The square roots
appear only in the bounds. In an accelerated call, take positive integers $u,v$ with
$\alpha=u/v$. Among the positive integers, $m$ is the largest satisfying
$m^v\le n^u$, and $k$ is the smallest satisfying $k^{144}\ge n^{25}$.
The dyadic scale for $\tau$ can likewise be selected by comparing
$\tau^{192}$ with $n^{23}L^{192}$.
The initial perturbation of a batch contributes one additional column,
so the low-rank subroutines use at most $k+1$ perturbation columns.

Apply Proposition~\ref{prop:seed} with tolerance $\delta$.
Denote the resulting matrix by $B_{\rm s}$ and
scale its right columns to lengths in $[1/2,1]$, applying the inverse
scalings to the corresponding dual rows. Its distance from the input is less than
$\delta/4$, and its minimum gap is at least $\delta/(8n)$.
Initialise $B=B_0=B_{\rm s}$ and $\mathcal E=0$. In addition to the
current matrix $B$, retain the matrix $\mathcal E$, which approximates
$B-B_{\rm s}$. We use the potential
\begin{equation}\label{eq:local-potential}
 \Phi(B,\mathcal E)=\log|D(B)|-\nu\|\mathcal E\|_F^2
\end{equation}
to control both the displacement and the number of updates.
To evaluate the phase-dependent part of the penalty during a batch,
we need only a Frobenius inner product. Let $B_0=X\Lambda X^{-1}$ be the reference matrix,
with $\Lambda=\operatorname{diag}(\lambda_1,\ldots,\lambda_n)$.
Let $x_i$ be column $i$ of $X$ and $w_i^*$ row $i$ of $X^{-1}$.
At each change of reference, update $B_0,X,\Lambda$. Between these
changes, express each exact rank-one update in reference coordinates.
The ratios $K_i$ refer to the matrix currently under consideration, whereas the
precision exponents $b_i$ remain attached to the fixed reference
eigensystem throughout a batch.
To change the reference, apply
Proposition~\ref{prop:batch-adoption}, requiring a correction of
Frobenius norm at most $\eta$ and a discriminant modulus at least
$0.9$ of its value before the correction. The matrix $\mathcal E$
is unchanged by this correction.

\subsection{The local algorithm}

The choice between an individual update and a batch depends on a sum
of weights computed in the reference basis. When this sum is large,
the increase in $\Phi$ bounds the total cost of recomputing the
eigensystem. When it is small, it bounds the precision needed to reuse
the reference basis.

In an accelerated call, repeat the following steps.
\begin{enumerate}
\item Normalise the right columns to lengths in $[1/2,1]$ and compute
\[
 K_i^2=\frac{\|x_i\|^2\|w_i\|^2}
              {\min_{j\ne i}|\lambda_i-\lambda_j|^2}
\]
exactly. Here $s_i=\min_{j\ne i}|\lambda_i-\lambda_j|$ is the nearest
gap in the reference spectrum. All threshold and weight comparisons
use squared quantities, so no square roots are required.
If at most $m$ values exceed $Q^2$, return the reference
matrix and its eigensystem. Otherwise, assign weights
\[
 \gamma_i=
 \begin{cases}
 \lfloor\log_2(rK_i/(2^{16}n))\rfloor,&K_i>2^{20}n/r,\\
 0,&K_i\le2^{20}n/r,
 \end{cases}
 \qquad W=\sum_i\gamma_i .
\]
For a positive weight, compare $(rK_i/(2^{16}n))^2$ with powers of
four. Successive squaring brackets the exponent, and the stored powers
in descending order determine the largest integer exponent whose power
does not exceed the squared ratio.
\item If $W\ge n$, perform the individual update below, change the
reference, and return to the first step.
\item If $W<n$, preprocess both reference bases as in
Section~\ref{sec:mixed-range-bridge}, including the tables for the sign
directions of Lemma~\ref{lem:batch-signs}. Apply the initial perturbation
of Lemma~\ref{lem:old-poles} with radius parameter $\min(r,\eta)$.
This perturbation has norm at most $\eta$ and is omitted from
$\mathcal E$. Perform at most $k$ updates in the batch. Before each
update, compute the characteristic polynomial and approximate its roots.
At each approximate eigenvalue, estimate the norms of the left and
right eigenvectors, their pairing and the distance to the other roots,
as in Section~\ref{sec:mixed-range-bridge}. The resulting estimates
satisfy $K_i/4\le\widehat K_i\le4K_i$ for the eigenvalues of the
current matrix. Select the indices with $\widehat K_i>4M$.
If at most $m$ indices are selected, retain all other indices,
change the reference subject to the tests below, and return its
eigensystem. Otherwise, reconstruct any $m$ selected pairs and perform
a batch update. After $k$ updates, change the reference and return
to the first step.
\end{enumerate}
A terminal call instead uses the procedure below until every
$K_i\le Q_{\rm term}$. The final subsection determines whether to
recurse and how to assemble the output.

\paragraph{Directions for an individual update.}
Choose the signs in the first assertion of Lemma~\ref{lem:signs}
with weights $\gamma_i$. Let $s$ be the least integer power of two
with $s^2\ge n$, so $\sqrt n\le s<2\sqrt n$, and rescale both
sign vectors to have entries $\pm1/(2s)$. Their norms are at most
$1/2$, and
\[
 |v^*P_i u|\ge \|P_i\|/(64n)
\]
on indices carrying at least $W/25$ total weight. In
$\mathbf a=X^{-1}u$ and $\mathbf b^*=v^*X$, replace each entry of modulus below
\[
 \theta=\bigl[100n^3(1+\|X\|_F^2+\|X^{-1}\|_F^2)\bigr]^{-1}
\]
by $\theta$, and reconstruct the two directions. The Euclidean
change in each direction is at most
\[
 2\sqrt n\,\theta\max(\|X\|_F,\|X^{-1}\|_F)
 \le \frac1{100n^{5/2}}.
\]
Thus both directions have norm at most one. On the indices selected by
the sign construction, each coordinate has modulus at least the
corresponding column or dual-row norm divided by $4s$. These coordinates
already exceed $\theta$, so the overlap estimate still holds. The nonzero coordinates
also make the new spectrum disjoint from the reference spectrum
at every nonzero parameter for which it is simple.

\paragraph{Directions for a batch update.}
For the $m$ selected pairs, Lemma~\ref{lem:batch-signs} supplies
nominal vectors $u_0,v_0$ with entries $\pm1/(2s)$ and exact
reference factors $\mathbf a,\mathbf b^*$. Define the actual directions by
$u=X\mathbf a$, $v^*=\mathbf b^*X^{-1}$. Both actual and nominal directions have
norm at most one, and on at least $m/25$ selected labels,
\[
 |v^*P_i u|\ge \|P_i\|/(2^{13}n).
\]
The errors in the original coordinates are at most
\[
 \epsilon=2^{-40}n^{-6}(1+L)^{-2}.
\]
The exact reference factors suffice during the batch.
For an accepted parameter $t$, update $B$ by $tuv^*$ and
$\mathcal E$ by $tu_0v_0^*$. The parameter bound below gives
\[
 \|tuv^*-tu_0v_0^*\|_F
 \le 3r\epsilon\le\eta.
\]
During a batch, $\mathcal E$ is represented by its explicit value at the
start and the subsequent nominal rank-one terms. At the end of the
batch, form the explicit matrix as in Lemma~\ref{lem:proxy-query}.
If the batch ends early, also require $K_i^2\le Q^2$ on the corrected
columns corresponding to the retained labels.

\paragraph{The phase choice.}
For a proposed update, let $F=uv^*$ denote its actual direction
matrix. Denote the direction used to update $\mathcal E$ by $F_0$.
Thus $F_0=F$ for an individual or terminal update, and
$F_0=u_0v_0^*$ in a batch. The cross term in the quadratic penalty
determines the coefficient $c$ in the phase search. Define
\[
 g(z)=D(B+rzF),\qquad
 c=2\nu r\,\operatorname{tr}(\mathcal E^*F_0).
\]
For an individual or terminal update compute $c$ directly and set
$c_0=c$. In a batch, apply Lemma~\ref{lem:proxy-query} with real
multiplier $2\nu r$ to obtain $c_0$ with $|c_0-c|\le1/100$.
This calculation uses the explicit matrix at the start of the batch and
the subsequent nominal updates. The proof of Proposition~\ref{prop:stage}
verifies the hypotheses of the lemma.
Apply Proposition~\ref{prop:weighted-phase} to $g$ and $c_0$.
For a batch update, use the reference nodes $\lambda_i$ as its additional
nodes. The initial perturbation and the previous phase choices
ensure that the current characteristic polynomial $p$ satisfies
$p(\lambda_i)\ne0$. No additional nodes are needed for individual
or terminal updates, whose adjusted directions already ensure
disjointness from the reference spectrum. The selected phase retains
at least $0.9$ of the maximum of
$|g(z)e^{-cz}|$ on the unit circle, with the following radius bound:
\[
 1\le |z|\le1+\frac1{500J},\qquad
 J=2n+64A,
\]
where $A$ is the least power of two not smaller than
$1+|\operatorname{Re}c_0|+|\operatorname{Im}c_0|$.
Use $t=rz$. Thus $|t|<2r$. For an individual update add $tF$ to $\mathcal E$.
The Taylor polynomial in the proposition allows the algorithm to select
the phase using field operations and comparisons.

\paragraph{Terminal updates.}
At a reference matrix choose any index with $K_i>Q_{\rm term}$,
if one exists. Compute upper norm estimates $c_w,c_x$ within a
factor $1.01$ for the norms of its dual row and right column, and take
\[
 u=\frac{w_i}{2c_w},\qquad v=\frac{x_i}{2c_x}.
\]
Both norms are at most $1/2$, and
$|v^*P_i u|\ge\|P_i\|/8$. Apply the same coordinate adjustment
as for an individual update. The two targeted coordinates exceed
$\theta$, so the overlap is preserved and both adjusted norms
are at most one. Use the weighted phase choice with $F_0=F$,
update $\mathcal E$ by the actual rank-one change, and change
the reference. If no such index exists, return the entire
eigensystem.

\subsection{Gain, precision and termination}

At the start of a batch we know that $W<n$. This bounds the sum of the
precisions needed for the individual reference eigenvectors.

\begin{lemma}[Precision for the reference eigensystems]\label{lem:total-precision}
At the beginning of a batch, the column exponents of
Section~\ref{sec:mixed-range-bridge} may be chosen with
$\sum_i b_i=O(nL)$ simultaneously for both reference bases.
The accuracy $\epsilon$ required for the sign directions gives the same
total preprocessing bound.
\end{lemma}
\begin{proof}
Choose $b_i$ to be the first power of two above
$20\lceil\log_2(2n)\rceil$ plus upper integer estimates of the
two nonnegative logarithms of the dual norm and reciprocal gap.
For the dual norm use the larger requirement for the separately
scaled right and left bases. Their primal column lengths are in
$[1/2,1]$, so their dual norms lie between $\|P_i\|$ and
$2\|P_i\|$. Since $\|P_i\|\ge1$ and $s_i\le4$, the two
nonnegative logarithms differ in sum from $\log_2^+K_i$ by at
most an absolute constant. The definition of $r$ gives
$\log_2(1/r)=O(L)$. Below the weight threshold this yields
$O(L)$, and above it the upper bound is $O(L)+\gamma_i$.
Dyadic rounding costs a factor of two. Consequently
$\sum_i b_i=O(nL+W)=O(nL)$.
The tables for the sign directions require only an additional
$O(\log n+\log(1+L))=O(L)$ bits per row or column, which
preserves this sum bound.
\end{proof}

\begin{lemma}[Gain per update]\label{lem:gain}
An individual update increases $\Phi(B,\mathcal E)$ by at least
$W/100$, including its matrix correction. Each accepted batch update
before the correction gains at least $m/25$. Each terminal update
gains at least one, including its correction. The total loss from
the initial perturbation and final correction of a batch is less
than two whenever the reference discriminant modulus is at least $G_0$.
\end{lemma}
\begin{proof}
The mean estimate in \eqref{eq:ascent}, converted to natural
logarithms, and $\int_0^{2\pi}\operatorname{Re}(ce^{it})\,dt=0$
show that the circle mean of $\log|g(z)e^{-cz}|-\log|g(0)|$
has the same lower bound as the mean discriminant gain.
The loss in logarithmic modulus from the true weighted maximum is
at most $-\log(0.9)<1/4$, by Proposition~\ref{prop:weighted-phase}. Moreover,
\[
 \nu |t|^2\|F_0\|_F^2\le 2\nu r^2\le 2^{-15}m.
\]
Here the sharper bound $|t|/r\le1+1/(500J)$ is used in the
first inequality. Subtracting this quadratic term gives the
increment of $\Phi$ exactly.

For an individual update, the indices on which both overlap bounds
hold have total weight at least $W/25$. At each such index of positive weight,
\[
 \frac{r|v^*P_i u|}{8s_i}
 \ge\frac{rK_i}{512n}\ge2^{\gamma_i}.
\]
The mean logarithmic gain is thus at least $(\log2)W/50$.
Since $m\le n\le W$ and $W\ge2^{20}$, the phase loss,
quadratic term, and final factor $0.9$ leave a gain at least
$W/100$. The same lower bound holds before the correction.

For a batch update, a selected true ratio exceeds $M$.
On at least $m/25$ such indices,
\[
 (r/16)|v^*P_i u|>4s_i.
\]
Lemma~\ref{lem:collision} gives at least $m/50$ zeros of $g$
in $|z|\le1/8$. Each contributes at least $3\log2$ to
the circle mean. Since $m\ge2^{15}$ in an accelerated
call, the phase and quadratic losses leave a gain at
least $m/25$.

For a terminal update, $K_i>2^{12}/r$ and the overlap
$\|P_i\|/8$ give the same strict collision inequality.
Integrality gives at least one zero of $g$ in $|z|\le1/8$.
Its contribution $3\log2$, less the phase and quadratic
losses and the correction loss, is greater than one.
Finally, the initial perturbation of a batch retains at least
$1-1/(500n^3)$ of its discriminant modulus, and the correction
retains $0.9$. Neither changes $\mathcal E$, so their total
loss is less than two.
\end{proof}

\begin{proposition}[Termination of a local stage]\label{prop:stage}
The matrix returned by the local algorithm has simple spectrum and lies
within $\delta/3$ of its input. In an accelerated call, at most $m$ indices have
$K_i>Q$. In a terminal call every $K_i\le Q_{\rm term}$.
Let $\rho=545/576$. The local algorithm costs
$\widetilde O(n^3+n^{3+\rho-\alpha})$ in the accelerated case
and $\widetilde O(n^4)$ in the terminal case, with a fixed power
of $1+\log(n/\delta)$ suppressed.
\end{proposition}
\begin{proof}
Consider any initial segment of at most $T$ matrix changes, counting initial
perturbations, accepted rank-one updates and corrections separately.
We first bound the distance from the input and the required precision
throughout this segment. At a change of reference, the norm and
discriminant estimates for the trial matrix determine the accuracy
needed to compute its correction. We then show that the algorithm
stops before $T$ changes.

\paragraph{Displacement bounds.}
By construction,
\begin{equation}\label{eq:proxy-error}
 \|B-B_{\rm s}-\mathcal E\|_F\le T\eta=\delta/2^{14}.
\end{equation}
At a batch update, Lemma~\ref{lem:batch-signs} bounds the difference
between the actual and nominal updates by $\eta$, including the factor
$|t|$. At an individual or terminal update the change to $\mathcal E$
is exact. Each initial perturbation and each correction is omitted
from $\mathcal E$ and has Frobenius norm at most $\eta$.

The initial matrix satisfies $|D(B_{\rm s})|\ge16G_0$.
Lemma~\ref{lem:gain} gives
\[
 \Phi(B,\mathcal E)\ge\log|D(B_{\rm s})|-2
\]
whenever we compute roots or an eigenbasis. Individual updates and
batches of the prescribed length increase the potential. Only the
last batch in the segment can be shorter, and its initial perturbation
and final correction decrease the potential by less than two in total.
As long as $\|B\|<2$, we have $|D(B)|\le4^d$, whence
\[
 \nu\|\mathcal E\|_F^2
 \le \log\frac{4^d}{|D(B_{\rm s})|}+2
 < 2n^2L+2.
\]
The definition of $\nu$ therefore gives
$\|\mathcal E\|_F<\delta/2^{10}$.
Together with \eqref{eq:proxy-error} and the initial perturbation bound
$\delta/4$, this bounds the distance from the input by
\[
 \delta/4+\delta/2^{10}+\delta/2^{14}<\delta/3.
\]
Every trial change has norm at most $2r<\delta/2^{18}$,
and the initial batch perturbations and corrections have norm at most $\eta$.
To show that the norm remains below two, suppose the current matrix
is within $\delta/3$ of the input. The next trial matrix has norm at most
$1+\delta/3+2r+\eta<2$. The phase estimate applies to this matrix
before its roots are approximated. After the
correction tests pass, the correction has Frobenius norm at most $\eta$
and preserves at least $0.9$ of the discriminant modulus. Applying the
potential estimate again gives the displacement bound after the correction.
The strict norm bound also places the matched approximants in
$|z|<2$ at the precision used below.

Since the quadratic term is nonnegative, the same potential
bound implies
$|D(B)|\ge e^{-2}|D(B_{\rm s})|>G_0$.
For the matrix at the selected phase, use the gain before correction.
Its norm and discriminant therefore give the bounds needed to approximate
its roots. Once the correction passes the discriminant test, we include
its additional loss of at most $-\log(0.9)$ in the potential estimate.

\paragraph{Precision bounds.}
Let $\mu_1,\ldots,\mu_n$ be the current eigenvalues.
Every current gap is at least $g=G_0/4^d$. The projector formula
\[
 P_i=\prod_{j\ne i}(B-\mu_jI)/(\mu_i-\mu_j)
\]
gives $\|P_i\|\le(4/g)^{n-1}$. For normalised reference
columns,
\[
 \log_2(1/g)\le4n^2L,\qquad
 \log_2(2+\|X\|+\|X^{-1}\|)\le8n^3L.
\]
At an initial perturbation, the characteristic polynomial
at each reference node has modulus
$|\sigma\alpha\beta p'_0(\lambda_i)|$, with
$|p'_0(\lambda_i)|\ge g^{n-1}$, in the notation of
Lemma~\ref{lem:old-poles}. Replacing its radius parameter by
$\min(r,\eta)$ increases the bound on $\log_2(1/\sigma)$ by at most
$O(L+\log(1+L))$. At an individual or terminal update, the adjusted
coordinates give the lower bound $r\theta^2|p'_0(\lambda_i)|$.
During a batch, each radial adjustment divides the
previous lower bound by at most $2^{16}nJ$.

At the start of every batch, $\|\mathcal E\|_F<\delta$, and
\[
 2\nu r\delta\le2^5 n\sqrt{mL}\le2^{10}n^2(1+L).
\]
Thus the real multiplier $2\nu r$ satisfies the hypotheses of
Lemma~\ref{lem:proxy-query}. The bound on $\|\mathcal E\|_F$ also gives
$|c|\le2\nu r\|\mathcal E\|_F\|F_0\|_F
=O(n\sqrt{mL})$, so $J=O(n\sqrt{mL}+n)$.
Consequently the logarithm of the reciprocal of the lower bound for
the characteristic polynomial values is at most
\[
 O\bigl(n^3L+T\log(n(1+L))\bigr).
\]
Since every such polynomial value is a product of $n$ distances, each
at most four, dividing its modulus by $4^{n-1}$ gives a lower bound
on the distance to every current eigenvalue. The middle factor in
the factorisation of $X^{-1}(B-B_0)X$ from
Section~\ref{sec:transfer-rebase} has norm at most
$4\|X^{-1}\|\|X\|\|B-B_0\|$.
To compare these estimates with $H$, observe that both
$|\sigma\alpha\beta p_0'(\lambda_i)|$ and $r\theta^2|p_0'(\lambda_i)|$
are at least $2^{-160n^3L}$. Indeed,
$|p_0'(\lambda_i)|\ge g^{n-1}$, each squared Frobenius norm in the
denominators is at most $n$ times the squared operator norm, and
$\log_2(1/\sigma)$ is bounded by the larger of the two exponents in
the definition of the initial batch perturbation. Moreover $J\le 2^{16}n^2L$,
so each radial adjustment increases the logarithm of the reciprocal
lower bound by at most $10L$. The logarithm of the reciprocal distance
is therefore bounded by
$160n^3L+10TL+2n<(256nL)^8$. The same estimates bound the norm of
the middle factor by $H$.

Matched approximants in the disc of radius two are obtained
by applying the scalar algorithm to $2^{-n}p(2z)$ at half
the required absolute error. For matched error $e$, a sufficient coefficient exponent for this
scaled polynomial is
$n(\lceil\log_2(2/e)\rceil+\lceil\log_2(4n)\rceil)$, by
Theorem~\ref{thm:paper-roots}. To obtain a correction of Frobenius norm at most $\eta$, it suffices
to take $\lceil\log_2(1/e)\rceil$ to be a sufficiently large absolute
multiple of $n\log H+\log(1/\eta)$.
Thus the logarithm of the coefficient exponent is $O(\log(nL))$.
The required powers of two are constructed by
$O(\log(nL))$ repeated squarings. The individual exponents used by the fixed
reference tables are those of Lemma~\ref{lem:total-precision}.
The tables for the Frobenius pairings require only
$O(\log n+\log(1+L))$ additional bits to approximate the phase
coefficient to absolute error at most $1/100$.

\paragraph{Tests at a change of reference.}
Consider first the discriminant test at a change of reference.
For gaps at least $H^{-1}$, moving each node by at most $e$
changes each pairwise distance by at most $2e$, a relative error
of at most $2He$. The inequality
$(1-2He)^d\ge1-2dHe$ shows that the corrected discriminant
retains at least $0.9$ of its previous modulus if $2dHe<1/10$.
When a batch ends early, we must also check the ratios for the
retained eigenvalues. Each has $K_i\le16M$ before the correction.
Keeping the same indices, scale each corresponding exact
column by its computed divisor, so its length is in
$[1/2,2]$. If its error is at most $e'$, the matrix error
is at most $\sqrt n e'$, while the true scaled matrix $X_t$ has inverse norm at most
$2\sqrt n(4H)^{n-1}$. The Neumann series for $(I+X_t^{-1}(X-X_t))^{-1}$ gives
\[
 \|X^{-1}-X_t^{-1}\|
 \le2\|X_t^{-1}\|^2\sqrt n e'
\]
provided $\|X_t^{-1}\|\sqrt n e'<1/2$.
Each true column and dual row has norm at least $1/2$. By the displayed
inverse estimate and Lemma~\ref{lem:transfer-vectors}, a sufficiently
large absolute exponent in the root accuracy makes the computed column
and dual-row norms at most $1.1$ times their true norms. Taking also
$2He\le1/11$ bounds the reciprocal nearest gap by $1.1$ times its old
value. Since $1.1^3<3/2$, every retained ratio is below $24M<Q$
after correction.
If refinement is needed, the labels are recovered from disjoint
radius-$1/(4H)$ discs about the first eigenvalue approximants used in
the resolvent formulas.
Test all pairs at the change of reference in $O(n^2)$ work.
Distinctness of the trial nodes and disjointness from the reference nodes
can be tested together by lexicographically sorting the two node
lists in $O(n\log n)$ comparisons before forming a Cauchy denominator. The correction norm is
tested from its thin factors. The old discriminant is available from the
current monic characteristic polynomial by a scalar resultant, and the
new one from the pairwise differences of the trial nodes. Squared moduli
suffice for the comparison. These tests cost at most
$\widetilde O(n^2(k+1)^{\omega-2}+nk^5)$ at a change of reference.
The fixed precision above suffices on the first attempt. Alternatively,
starting with the accuracy of Lemma~\ref{lem:transfer-vectors} and doubling the requested exponent
until the tests pass multiplies this cost by at most $O(\log(nL))$.

\paragraph{Termination.}
Let $q$ count individual, batch, and terminal updates
in the initial segment under consideration. A complete batch gains at least one per rank-one
update after its two losses, and an unfinished batch loses
at most two. If the segment ends between the phase choice and the
correction of an individual or terminal update, that update has
already increased the potential by at least one.
Hence
\[
 \Phi(B,\mathcal E)-\Phi(B_{\rm s},0)\ge q-2.
\]
On the other hand this increase is less than $2n^2L$.
Thus $q<2n^2L+2$. Every batch completed at its prescribed length contains $k$ updates.
Only the final batch can return earlier and contain none. Counting its
initial perturbation and final correction separately gives
at most $4q+4$ matrix changes, strictly fewer than $T$.
Thus the algorithm cannot reach $T$ matrix changes. Each inner procedure
terminates, so the local stage terminates. At either stopping test the
ratios satisfy the stated bounds, and the matrix lies within
$\delta/3$ of the input.

\paragraph{Operation count.}
Use $\omega<2.38$ from
\cite{alman-multiplication}, $a=23/192$, $b=25/144$ and
$\rho=545/576$. Let $m=\lfloor n^\alpha\rfloor$ and
$k=\lceil n^b\rceil$, where $\alpha$ is the exponent for the current
depth. The preparation and change of reference for a batch cost
\[
 \widetilde O\bigl(n^{2+a}+n^2k^{\omega-2}+nk^5\bigr).
\]
This includes the packed reference transformations and tables for the Frobenius pairings.
Updating the accumulated matrix $\mathcal E$ at the end of a batch costs
$\widetilde O(n^2k^{\omega-2})$ by multiplying the arrays of stored sign directions as in Lemma~\ref{lem:proxy-query}. The table gives the cost of one update within the batch.
\begin{center}
\begin{tabular}{ll}
\hline
Computation & Cost up to fixed logarithmic factors\\
\hline
Discriminant, weighted phase and roots
 & $n^{(1+\omega)/2}+n^{1.9}+n\sqrt m$\\
Resolvent cofactors & $nk^5$\\
High-precision coordinates & $n^{2-a}k^{\omega-2}+nk^{\omega-1}$\\
Other norm estimates & $nk^2n^{5a}$\\
Selected vectors and signs & $nm$\\
Reference coordinates and penalty pairings & $n+nk$\\
\hline
\end{tabular}
\end{center}
Forming the numerator factors for the norm estimates costs
$\widetilde O(nk^2+nk)$ and is covered by the displayed bounds.
The identities
\[
 2-a+\frac{19}{50}b=1+2b+5a=2+a-b=1+\rho
\]
show that each row except the one for selected vectors and signs is
$\widetilde O(n^{1+\rho})$, as is the batch preparation divided by $k$.
For the remaining terms use $5b<\rho$,
$2+(\omega-3)b<1+\rho$, $1+(\omega-1)b<1+\rho$ and
$1+\alpha/2\le3/2<1+\rho$.
Thus a complete batch costs
$\widetilde O(k(n^{1+\rho}+nm))$ and gains $\Omega(km)$.
The total available increase in the potential is $O(n^2L)$, so these
batches cost $\widetilde O(n^{3+\rho-\alpha}+n^3)$ in total.
Individual updates cost $\widetilde O(n^2)$ and gain $\Omega(n)$, giving
another $\widetilde O(n^3)$.

There is at most one incomplete batch. Its preparation and change of
reference cost $\widetilde O(n^{2+a})$, since
$2+(\omega-2)b<2+a$ and $1+5b<2+a$. Its accepted updates have the
same ratio of operation count to potential increase as those in a
complete batch. The final test costs $\widetilde O(n^{1+\rho}+nm)$,
which is below the cubic bound. The initial perturbation and the exact reference tests cost
$O(n^2)$ per batch and are included in its preparation.
Forming the returned diagonalisation is included in the
change-of-reference cost. Forming the local output matrix from its
diagonalisation costs $O(n^3)$ once. Consequently the local
accelerated cost is
\[
 \widetilde O\bigl(n^3+n^{3+\rho-\alpha}+n^{2+a}\bigr).
\]
A terminal update costs $\widetilde O(n^2)$: the two directions and
the Frobenius pairing use the stored basis and the explicit accumulated
matrix, the pencil has width one, and the phase search has $m=1$.
There are $O(n^2L)$ such updates, giving $\widetilde O(n^4)$ work.
The initial eigensystem and the final ratio test are included in these bounds.
\end{proof}

\subsection{Triangular assembly and uniform bounds}
At the end of a local stage, the right eigenvectors for the indices
satisfying the cutoff span an invariant subspace. Changing the
trailing block can nevertheless change the projectors for the retained
eigenvalues. We choose a shift for that block and bound the change
using the gaps in the spectrum before the recursion.

\begin{proof}[Proof of Theorem~\ref{thm:main}]
There are at most two accelerated stages, with
\[
 \alpha_0=\frac{947}{1000},\qquad \alpha_1=\frac{39}{50},
 \qquad \rho=\frac{545}{576},
\]
followed, when needed, by the terminal procedure. At depth $j=0,1$, the
number of indices exceeding the local bound is at most
$m=\lfloor n^{\alpha_j}\rfloor$. A call below the
fixed dimension threshold uses the terminal procedure immediately.
We use the local cutoffs $Q$ and $Q_{\rm term}$ of
\eqref{eq:local-thresholds} and their bounds in \eqref{eq:local-threshold-bounds}.

\paragraph{The invariant subspace and the trailing block.}
At a nonterminal call, let $B$ be the matrix returned by the local
procedure. It has simple spectrum and lies within $\delta/3$ of the input.
Take the invariant span of all right eigenvectors whose ratios satisfy
$K_i\le Q$. Its codimension $\ell$ is at most $m<n$.
If $\ell=0$, return $B$ and its stored eigensystem.
Otherwise, exact Gram--Schmidt projections, completed with standard
coordinate vectors and scaled as in Section~\ref{sec:foundations}, give
a basis $S$ of orthogonal columns with $\kappa_2(S)\le1.02$ and
\[
 S^{-1}BS=\begin{pmatrix}G&C\\0&H_b\end{pmatrix}.
\]
Let $p=n-\ell$ be the order of $G$. Both diagonal blocks have norm at
most $1.02(1+\delta/3)<2$, so $H_b/2$ is a contraction.
Apply the algorithm recursively to $H_b/2$ with tolerance $\delta/64$
and multiply the returned matrix by two to obtain $H_{\rm rec}$.
Induction on the remaining depth gives
$\|H_{\rm rec}-H_b\|<\delta/32$.
For a scalar trailing block the recursive call leaves the entry
unchanged and gives the unit column as an eigenbasis.

\paragraph{A shift adapted to the leading gaps.}
Let $\lambda_1,\ldots,\lambda_p$ be the eigenvalues retained in $G$,
and let $\sigma_i$ be the nearest-eigenvalue gap of $\lambda_i$ in the
full old spectrum of $B$. Thus the original spectral projectors
$P_i^{\rm old}$ satisfy
\[
 \|P_i^{\rm old}\|\le Q\sigma_i.
\]
The discs $D(\lambda_i,\sigma_i/2)$ have disjoint interiors.
Since $\|B\|<2$ and $\sigma_i<4$, they lie in $D(0,4)$. Consequently
\begin{equation}\label{eq:assembly-gap-packing}
 \sum_{i=1}^p\sigma_i^2\le64,\qquad
 \sum_{i=1}^p\sigma_i\le8\sqrt p\le8\sqrt n .
\end{equation}

Let $\mu_1,\ldots,\mu_\ell$ denote the eigenvalues of $H_{\rm rec}$.
We choose a small complex shift $z$ so that
\begin{equation}\label{eq:assembly-adaptive-gap}
 |\lambda_i-(\mu_j+z)|
 \ge \frac{\delta}{2^{20}}
       \max\left\{\frac1n,\frac{\sigma_i}{\sqrt\ell}\right\}
 \quad(1\le i\le p,\ 1\le j\le\ell).
\end{equation}
For each $i$, choose $\lambda_i'$ minimising $|\lambda_i-\lambda|^2$
over $\lambda\in\operatorname{spec}(B)\setminus\{\lambda_i\}$ by field comparisons.
Set $\widetilde\sigma_i=|\Re(\lambda_i-\lambda_i')|+|\Im(\lambda_i-\lambda_i')|$. Then
$\sigma_i\le\widetilde\sigma_i\le2\sigma_i$.
Find an integer $t$ with $\sqrt\ell\le t\le2\sqrt\ell$ by comparisons
of integer squares. Define
\[
 \gamma=\frac{\delta}{2^{20}n},\qquad q=2^8n,\qquad
 r_i=\frac{\delta}{2^{20}}
       \max\left\{\frac1n,\frac{2\widetilde\sigma_i}{t}\right\}.
\]
Consider the shifts $z=\gamma(a+\mathrm i b)$ with integer
$-q\le a,b\le q$. For each pair $(i,j)$, mark the grid points in the
closed square of radius $r_i$ in the maximum norm about
$\lambda_i-\mu_j$.
For $R_i=r_i/\gamma\ge1$, this square contains at most
$25R_i^2$ grid points, and
\[
 R_i^2\le1+\frac{16n^2\sigma_i^2}{\ell}.
\]
By \eqref{eq:assembly-gap-packing}, the total number marked, counting
repetitions, is at most
\[
 25\left(p\ell+16n^2\sum_i\sigma_i^2\right)
 \le25(1/4+1024)n^2<2^{15}n^2 .
\]
There are $(2q+1)^2>2^{18}n^2$ candidate shifts. An unmarked point
therefore exists and satisfies \eqref{eq:assembly-adaptive-gap}.
Its norm is at most
$\sqrt2\,q\gamma=\sqrt2\,\delta/2^{12}<\delta/1024$.
The integer endpoints of each marked square are found by comparisons
within $[-q,q]$. Enumerating the marked points, sorting them, and finding
a missing point cost $O(n^2\log n)$ operations and comparisons.
This includes the $O(n^2)$ comparisons needed for the old nearest gaps.

Set $H_b'=H_{\rm rec}+zI$ and $\Delta=H_b'-H_b$. Then
\begin{equation}\label{eq:assembly-trailing-error}
 \|\Delta\|<\frac{33\delta}{1024}.
\end{equation}
In particular
$\|H_b'\|\le1.02(1+\delta/3)+33\delta/1024<2$.
The two diagonal spectra are simple and disjoint. Define
\[
 T_{\rm new}=\begin{pmatrix}G&C\\0&H_b'\end{pmatrix},\qquad B_{\rm out}=S T_{\rm new}S^{-1}.
\]
The distance from $B_{\rm out}$ to the input is less than
\[
 \delta/3+1.02(\delta/32+\delta/1024)<\delta.
\]
The leading eigenbasis is inherited from the selected eigenvectors, and
the recursive call supplies the trailing eigenbasis. In these bases,
the Sylvester equation is solved by dividing each entry by the
corresponding difference of a leading and a trailing eigenvalue.
This gives an exact full eigenbasis and its inverse
in $O(n^3)$ operations. The spectral tests and grid selection also use real comparisons.

\paragraph{The sum of the projector norms.}
For a matrix $A$ of positive order with simple spectrum, define
$\mathcal S(A)=\sum_{\lambda\in\operatorname{spec}(A)}\|P_\lambda^A\|$.
In particular a scalar has $\mathcal S=1$.
We use
\begin{equation}\label{eq:assembly-sum-conditioning}
 \kappa_V(A)\le\mathcal S(A).
\end{equation}
Indeed, choose unit right columns $x_i$ and dual rows $w_i^*$.
If $p_i=\|P_i^A\|=\|w_i\|$, multiplying column $i$ by $\sqrt{p_i}$
and its dual row by $1/\sqrt{p_i}$ makes both squared Frobenius norms
equal to $\sum_i p_i$. The operator norms are bounded by the Frobenius
norms. This rescaling is used only for the condition estimate.

For each retained eigenvalue, the old projector in the coordinates of
$S$ has the form
\[
 \Pi_i=S^{-1}P_i^{\rm old}S
       =\begin{pmatrix}U_i&R_i\\0&0\end{pmatrix},
 \qquad
 \|\Pi_i\|\le A_0\sigma_i,\qquad A_0=1.02Q,
\]
where $U_i$ is the corresponding projector of $G$.
Let $Y=\sum_iR_i\in\mathbb C^{p\times\ell}$ and $S_g=\sum_i\|\Pi_i\|$.
The projector identities give
\[
 \sum_iU_i=I_p,\qquad U_iY=R_i,\qquad
 C=GY-YH_b,\qquad
 \|Y\|\le S_g\le8A_0\sqrt n.
\]
These identities follow by summing the commutation relations for $\Pi_i$ and $S^{-1}BS$ and using $\Pi_i\Pi_j=0$ for $i\ne j$, $\Pi_i^2=\Pi_i$.

For an eigenvalue $\mu$ of $H_b'$, let $T_\mu$ be its spectral projector
and set $\mathcal S_{\rm ch}=\sum_{\mu\in\operatorname{spec}(H_b')}\|T_\mu\|$.
Doubling and shifting leave projectors unchanged, so this is the sum
returned by the recursive call. It is at least one.
Since $U_i C=R_i(\lambda_i I-H_b)$, the triangular resolvent gives
the following upper-right block for the new leading projector:
\[
 R_i'=R_i\bigl[I+\Delta(\lambda_iI-H_b')^{-1}\bigr].
\]
The upper-right block of the new trailing projector is
\[
 \begin{aligned}
  (\mu I-G)^{-1}CT_\mu
   &=-YT_\mu+(\mu I-G)^{-1}Y\Delta T_\mu,\\
  (\mu I-G)^{-1}Y
   &=\sum_i\frac{R_i}{\mu-\lambda_i}.
 \end{aligned}
\]
Define $F(\mu)=\sum_i\sigma_i/|\mu-\lambda_i|$. Expanding
$(\lambda_i I-H_b')^{-1}$ in the projectors $T_\mu$, the sum of the
norms of the corrections to the leading projectors is at most
$\|\Delta\|A_0\sum_\mu F(\mu)\|T_\mu\|$. The terms containing
$\Delta$ in the trailing projectors satisfy the same bound.
The remaining upper-right terms of the trailing projectors contribute at most
$\|Y\|\mathcal S_{\rm ch}$. Including the diagonal blocks gives
\[
 \mathcal S(T_{\rm new})
 \le S_g+
 \left(1+S_g+2\|\Delta\|A_0\sup_\mu F(\mu)\right)
       \mathcal S_{\rm ch}.
\]

To estimate $F$, choose a closest leading eigenvalue $\lambda_*$ to
$\mu$. Equation~\eqref{eq:assembly-adaptive-gap} bounds its contribution
by $2^{20}\sqrt\ell/\delta$.
For every other index,
$\sigma_i\le|\lambda_i-\lambda_*|\le2|\mu-\lambda_i|$.
Group these indices into shells
$r\le|\mu-\lambda_i|<2r$, with dyadic radii starting at $\gamma$.
For one such shell the disjoint discs
$D(\lambda_i,\sigma_i/4)$ lie in $D(\mu,3r)$, so
$\sum_{\rm shell}\sigma_i^2\le144r^2$.
If the shell contains $h$ indices, Cauchy--Schwarz bounds its contribution
to $F$ by $r^{-1}(h\sum_{\rm shell}\sigma_i^2)^{1/2}\le12\sqrt h$.
Both diagonal blocks after the trailing perturbation have norm below
two. Hence all these distances are less than eight, and at most
$L+17\le3L$ shells suffice. Cauchy--Schwarz across the shells yields
\begin{equation}\label{eq:assembly-weighted-resolvent}
 \sup_\mu F(\mu)
 \le\frac{2^{20}\sqrt\ell}{\delta}+21\sqrt{nL}.
\end{equation}
Using \eqref{eq:assembly-trailing-error}, $\ell\le n$, $\delta<1/2$,
and $S_g\le8A_0\sqrt n$, then conjugating back by $S$, gives
\begin{equation}\label{eq:assembly-sum-recurrence}
 \mathcal S(B_{\rm out})
 \le2^{18}Q\sqrt{nL}\,\mathcal S_{\rm ch}.
\end{equation}
The constant includes both factors of $1.02$.
If the trailing block is empty, the direct bound is
$\mathcal S(B)\le8Q\sqrt n$ instead.
No gap or condition number is assigned to an empty block.

\paragraph{Uniform conditioning and gap bounds.}
For the input dimension $N$ and tolerance $\delta_*$, define
\[
 L_0=\lceil\log_2(64N/\delta_*)\rceil .
\]
Let $n_j,\delta_j,L_j$ denote the dimension, tolerance and logarithmic parameter at depth $j$. Thus $\delta_j=\delta_*/64^j$ and, for $j\le2$,
$L_j\le L_0+12\le3L_0$.
The successive dimensions satisfy
\[
 n_0=N,\qquad
 n_1\le N^{947/1000},\qquad
 n_2\le N^{36933/50000}.
\]
In a terminal block, \eqref{eq:assembly-gap-packing} applied to its full
spectrum and the cutoff $Q_{\rm term}$ give
$\mathcal S\le8Q_{\rm term}\sqrt n$.
Together with \eqref{eq:local-threshold-bounds},
\eqref{eq:assembly-sum-recurrence}, and
$\lfloor n^\alpha\rfloor\ge n^\alpha/2$, this yields
\begin{equation}\label{eq:assembly-global-sum}
 \mathcal S(B_{\rm out})
 \le C_1\delta_*^{-3}L_0^{5/2}N^{256633/50000}.
\end{equation}
The three dimension contributions to this exponent are
\[
 \frac{5-\alpha_0}{2}=\frac{4053}{2000},\qquad
 \alpha_0\frac{5-\alpha_1}{2}=\frac{199817}{100000},\qquad
 \frac32\alpha_0\alpha_1=\frac{110799}{100000}.
\]
Their sum is $256633/50000$. Each accelerated factor
$Q\sqrt{nL}$ contributes one factor $\delta_j^{-1}L_j$, and the
terminal factor $Q_{\rm term}\sqrt n$ contributes
$\delta_j^{-1}L_j^{1/2}$. An empty trailing block removes subsequent
factors from the estimate. A scalar contributes one. If a call uses the
terminal procedure early because its dimension is below the fixed
threshold, its dimension-dependent factor is absorbed into the absolute
constant. These cases also satisfy \eqref{eq:assembly-global-sum}.

For a retained leading eigenvalue, $\|P_i^{\rm old}\|\ge1$ implies
$\sigma_i\ge Q^{-1}$. Discarding other eigenvalues can only increase its distance to the
nearest remaining neighbour when one exists. The corresponding terminal bound is
$1/\operatorname{gap}(B)\le Q_{\rm term}$ whenever the block has at least
two eigenvalues. A singleton has no internal pair.
Pairs from distinct diagonal blocks have separation at least $\delta/(2^{20}n)$.
Rescaling a recursive result in its parent call multiplies its internal
gaps by two. Subsequent shifts and changes of basis preserve them.
The powers of $N$ in the two accelerated bounds for the inverse gap are
\[
 2-\frac{\alpha_0}{2}=\frac{3053}{2000},\qquad
 \alpha_0\left(2-\frac{\alpha_1}{2}\right)
      =\frac{152467}{100000}.
\]
The terminal bound for the inverse gap has power $36933/50000$.
The bounds for pairs from distinct diagonal blocks have power at most
one. The first displayed power is the largest. Therefore
\begin{equation}\label{eq:assembly-global-gap}
 \operatorname{gap}(B_{\rm out})^{-1}
 \le C_2\delta_*^{-1}L_0^{1/2}N^{3053/2000}.
\end{equation}

Combining \eqref{eq:assembly-sum-conditioning},
\eqref{eq:assembly-global-sum}, and \eqref{eq:assembly-global-gap},
we obtain
\[
 \frac{\kappa_V(B_{\rm out})}{\operatorname{gap}(B_{\rm out})}
 \le C_3N^{166479/25000}\delta_*^{-4}L_0^3.
\]
The difference between $20/3$ and $166479/25000$ is
$\vartheta=563/75000>0$.
For $R=N/\delta_*>4$, we have $L_0\le7+\log_2R$, and
$(7+\log_2R)^3\le C_\vartheta R^\vartheta$.
Here $C_\vartheta$ is finite because
$(7+x)^3\,2^{-\vartheta x}$ is bounded for $x\ge2$.
Since $166479/25000>4$ and $\delta_*<1$, the preceding estimate is at
most $C R^{20/3}$, with an absolute constant $C$.

\paragraph{Operation count.}
There is at most one recursive call at each accelerated stage.
The local cost in dimension $n$ with exponent $\alpha_j$ is
$\widetilde O(n^{3+\max(0,\rho-\alpha_j)})$.
The top-level exponent is three because $\alpha_0>\rho$.
At depth one, its exponent as a power of $N$ is
\[
 \frac{947}{1000}
 \left(3+\frac{545}{576}-\frac{39}{50}\right)
 =\frac{43176571}{14400000}<3.
\]
The quartic terminal computation has exponent
$4\alpha_0\alpha_1=36933/12500<3$.
An earlier terminal call of bounded dimension costs at most a fixed
constant times its cubic cost.
The orthogonal bases and block assembly cost $O(n^3)$, and the adaptive
shift costs $O(n^2\log n)$.
All logarithmic parameters are $O(1+\log(N/\delta_*))$, and there are
at most three nontrivial calls. The total cost is therefore
$O(N^3\log^{c_{\log}}(N/\delta_*))$ for an absolute constant $c_{\log}$.
The perturbation bound holds at each return from the recursion.
Subtracting the input matrix from the final matrix gives $E$.
\end{proof}

\small

\bibliographystyle{ACM-Reference-Format}
\bibliography{deterministic-spectral-regularization}

\end{document}